\documentclass[a4paper]{article}

\usepackage[titletoc,title]{appendix}
\usepackage[cp1250]{inputenc}
\usepackage[IL2]{fontenc}
\usepackage{a4wide}
\usepackage[english]{babel}
\usepackage{euscript}
\usepackage{amstext,amsbsy,amscd,amssymb}
\usepackage{amsmath,enumerate}
\usepackage{amsfonts}
\usepackage{microtype}
\usepackage{color}
\include{diagxy}
\allowdisplaybreaks

\let\rarr=\rightarrow

\let\veps=\varepsilon
\let\mcal=\mathcal
\let\mfrak=\mathfrak
\let\eus=\EuScript

\def\N{\mathbb{N}}
\def\Z{\mathbb{Z}}

\def\C{\mathbb{C}}

\def\Hom{\mathop {\rm Hom} \nolimits}

\def\ad{\mathop {\rm ad} \nolimits}
\def\gr{\mathop {\rm gr} \nolimits}

\def\SL{\mathop {\rm SL} \nolimits}

\def\diag{\mathop {\rm diag} \nolimits}
\def\id{\mathop {\rm id} \nolimits}

\def\rank{\mathop {\rm rank} \nolimits}
\def\Ad{\mathop {\rm Ad} \nolimits}

\def\Ind{\mathop {\rm Ind} \nolimits}

\def\htt{\mathop {\rm ht} \nolimits}

\def\Cas{\mathop {\rm Cas} \nolimits}

\def\Specm{\mathop {\rm Specm} \nolimits}

\newsymbol\squares 1003

\long\def\proof #1{\noindent \emph{Proof.}\ #1 \hfill $\squares$
\medskip}

\newcounter{num}[section]
\numberwithin{equation}{section}
\numberwithin{num}{section}

\long\def\definition #1 {\refstepcounter{num} \noindent {\bf Definition \thenum.} #1

\medskip}

\long\def\theorem #1{\refstepcounter{num} \noindent {\bf Theorem \thenum.} #1

\medskip}

\long\def\lemma #1{\refstepcounter{num}  \noindent {\bf Lemma \thenum.} #1

\medskip}

\long\def\conjecture #1{\noindent {\bf Conjecture.}\ #1

\medskip}

\long\def\proposition #1{\refstepcounter{num}  \noindent {\bf Proposition \thenum.} #1

\medskip}

\long\def\corollary #1{\refstepcounter{num}  \noindent {\bf Corollary \thenum.}\ #1

\medskip}

\long\def\remark #1{\noindent {\bf Remark.}\ #1

\medskip}

\long\def\example #1{\noindent {\bf Example.}\ #1

\medskip}

\newenvironment{enum}{\begin{list}{}{\topsep=2pt \itemsep=0pt
\parsep=0pt}}{\end{list}}

\makeatletter
\newcommand*\if@single[3]{%
  \setbox0\hbox{${\mathaccent"0362{#1}}^H$}%
  \setbox2\hbox{${\mathaccent"0362{\kern0pt#1}}^H$}%
  \ifdim\ht0=\ht2 #3\else #2\fi
  }
\newcommand*\rel@kern[1]{\kern#1\dimexpr\macc@kerna}
\newcommand*\widebar[1]{\@ifnextchar^{{\wide@bar{#1}{0}}}{\wide@bar{#1}{1}}}
\newcommand*\wide@bar[2]{\if@single{#1}{\wide@bar@{#1}{#2}{1}}{\wide@bar@{#1}{#2}{2}}}
\newcommand*\wide@bar@[3]{%
  \begingroup
  \def\mathaccent##1##2{%
    \if#32 \let\macc@nucleus\first@char \fi
    \setbox\z@\hbox{$\macc@style{\macc@nucleus}_{}$}%
    \setbox\tw@\hbox{$\macc@style{\macc@nucleus}{}_{}$}%
    \dimen@\wd\tw@
    \advance\dimen@-\wd\z@
    \divide\dimen@ 3
    \@tempdima\wd\tw@
    \advance\@tempdima-\scriptspace
    \divide\@tempdima 10
    \advance\dimen@-\@tempdima
    \ifdim\dimen@>\z@ \dimen@0pt\fi
    \rel@kern{0.6}\kern-\dimen@
    \if#31
      \overline{\rel@kern{-0.6}\kern\dimen@\macc@nucleus\rel@kern{0.4}\kern\dimen@}%
      \advance\dimen@0.4\dimexpr\macc@kerna
      \let\final@kern#2%
      \ifdim\dimen@<\z@ \let\final@kern1\fi
      \if\final@kern1 \kern-\dimen@\fi
    \else
      \overline{\rel@kern{-0.6}\kern\dimen@#1}%
    \fi
  }%
  \macc@depth\@ne
  \let\math@bgroup\@empty \let\math@egroup\macc@set@skewchar
  \mathsurround\z@ \frozen@everymath{\mathgroup\macc@group\relax}%
  \macc@set@skewchar\relax
  \let\mathaccentV\macc@nested@a
  \if#31
    \macc@nested@a\relax111{#1}%
  \else
    \def\gobble@till@marker##1\endmarker{}%
    \futurelet\first@char\gobble@till@marker#1\endmarker
    \ifcat\noexpand\first@char A\else
      \def\first@char{}%
    \fi
    \macc@nested@a\relax111{\first@char}%
  \fi
  \endgroup
}
\makeatother

\newcommand\rsmraise[1]{%
  \ifx#1\displaystyle .8\else
    \ifx#1\textstyle .8\else
      \ifx#1\scriptstyle .6\else
        .45%
      \fi
    \fi
  \fi}

\title{Geometric construction of Gelfand--Tsetlin modules \\ over simple Lie algebras}

\author{Vyacheslav Futorny, Libor Křižka}

\AtEndDocument{\bigskip{\footnotesize%
  (V.\,Futorny) \textsc{Instituto de Matemática e Estatística, Universidade de Sa\~{o} Paulo, Caixa Postal 66281,\\ Sa\~{o} Paulo, CEP 05315-970, Brasil} \par
  \textit{E-mail address}: \texttt{futorny@ime.usp.br} \par
  \addvspace{\medskipamount}
  (L.\,Křižka) \textsc{Instituto de Matemática e Estatística, Universidade de Sa\~{o} Paulo, Caixa Postal 66281,\\ Sa\~{o} Paulo, CEP 05315-970, Brasil} \par
  \textit{E-mail address}: \texttt{krizka.libor@gmail.com} \par
}}

\date{}

\begin{document}

\maketitle

\begin{abstract}

In the present paper we describe a new class of Gelfand--Tsetlin modules for an arbitrary complex simple finite-dimensional Lie algebra $\mathfrak{g}$ and give their geometric realization as the space of `$\delta$-functions' on the flag manifold $G/B$ supported at the $1$-dimensional submanifold. When $\mathfrak{g}=\mathfrak{sl}(n)$ (or $\mathfrak{gl}(n)$) these modules form a subclass of  Gelfand-Tsetlin modules with infinite dimensional weight subspaces.
We  discuss their properties and describe the simplicity criterion for these modules in the case of the Lie algebra $\mathfrak{sl}(3,\mathbb{C})$.



\end{abstract}

\thispagestyle{empty}

\tableofcontents


\section*{Introduction}
\addcontentsline{toc}{section}{Introduction}

A study of simple \emph{weight} modules for a complex reductive finite-dimensional Lie algebra $\mfrak{g}$ is a classical problem in representation theory. Such modules have the diagonalizable action of a certain Cartan subalgebra $\mfrak{h}$ of $\mfrak{g}$. However, a complete classification of simple weight modules is only known for the Lie algebra $\mathfrak{sl}(2,\C)$ when the result is obvious. A classification of simple $\mathfrak{sl}(2,\C)$-modules was obtained in \cite{Block1981} and a classification of simple weight modules with finite-dimensional weight spaces for any $\mathfrak{g}$ was described in \cite{Fernando1990} and \cite{Mathieu2000}. A larger class of weight modules containing all modules with finite-dimensional weight spaces was extensively studied by many authors during several decades for the Lie algebra $\mfrak{g}$ of type $A$. The corresponding modules are called \emph{Gelfand--Tsetlin modules}. They are defined in a similar manner as the weight modules but instead of a Cartan subalgebra $\mfrak{h}$ of $\mfrak{g}$ one uses a certain maximal commutative subalgebra $\Gamma$ of the universal enveloping algebra $U(\mfrak{g})$, called the \emph{Gelfand--Tsetlin subalgebra}, see \cite{Drozd-Futorny-Ovsienko1994}.
Let us recall that Gelfand--Tsetlin modules are connected with the study of Gelfand--Tsetlin integrable systems \cite{Guillemin-Sternberg1983}, \cite{Kostant-Wallach2006}, \cite{Graev2007},
\cite{Colarusso-Evens2014}. Properties and explicit constructions of Gelfand--Tsetlin modules for the Lie algebra $\mfrak{g}$ of type $A$ was studied in \cite{Drozd-Ovsienko-Futorny1991}, \cite{Drozd-Futorny-Ovsienko1994}, \cite{Mazorchuk2001}, \cite{Ovsienko2002},  \cite{Graev2004}, \cite{Futorny-Grantcharov-Ramirez2014}, \cite{Futorny-Grantcharov-Ramirez2015}, \cite{Futorny-Grantcharov-Ramirez2016}, \cite{Futorny-Grantcharov-Ramirez2016a}, \cite{Futorny-Ramirez-Zhang2016}, \cite{Zadunaisky2017}, \cite{Vishnyakova2017a}, \cite{Vishnyakova2017b}, \cite{Ramirez-Zadunaisky2017} among the others.

Allowing only some generators of the Gelfand--Tsetlin subalgebra to have torsion on simple $\mathfrak{gl}(n,\C)$-modules leads to \emph{partial Gelfand--Tsetlin modules} for $\mathfrak{gl}(n,\C)$ which were studied in \cite{Futorny-Ovsienko-Saorin2011}. These modules decompose into the direct sum of $\Gamma$-submodules parameterized by prime ideals of the Gelfand--Tsetlin subalgebra $\Gamma$.

In the present paper, we study a special class of Gelfand--Tsetlin modules for $\mathfrak{gl}(n)$ and their analogues for an arbitrary complex simple finite-dimensional Lie algebra $\mathfrak{g}$. Let $\mfrak{h}$ be a Cartan subalgebra of $\mfrak{g}$ and $\mfrak{b}$ a Borel subalgebra of $\mfrak{g}$ containing $\mfrak{h}$. The Beilinson--Bernstein correspondence provides a geometric realization of Verma modules $M^\mfrak{g}_\mfrak{b}(\lambda)$ for $\lambda \in \mfrak{h}^*$. By this correspondence, to the Verma module $M^\mfrak{g}_\mfrak{b}(\lambda)$ we associate the vector space of `$\delta$-functions' on the flag manifold $G/B$ supported at the point $eB$. From this point of view, we introduce a family of  Gelfand--Tsetlin modules $W^\mfrak{g}_\mfrak{b}(\lambda)$ for $\lambda \in \mfrak{h}^*$ as the vector space of `$\delta$-functions' on the flag manifold $G/B$ supported at the $1$-dimensional submanifold going through the point $eB$. Modules  $W^\mfrak{g}_\mfrak{b}(\lambda)$ provide
first examples of Gelfand-Tsetlin modules for an arbitrary simple finite-dimensional $\mfrak{g}$ with respect to certain commutative subalgebra $\Gamma$ of the universal enveloping algebra $U(\mathfrak{g})$.

We will not discuss this geometric approach in details, instead we  describe everything from the algebraic point of view. However, geometry serves as a motivation for the purely algebraic construction of $\mfrak{g}$-modules $W^\mfrak{g}_\mfrak{b}(\lambda)$ considered in this paper.

\begin{table}[ht]
\centering
\renewcommand{\arraystretch}{1.5}
\begin{tabular}{c|c}
  Verma module & Gelfand--Tsetlin module \\
  \hline
  $M^\mfrak{g}_\mfrak{b}(\lambda)$ & $W^\mfrak{g}_\mfrak{b}(\lambda)$ \\
  weights $\lambda-Q_+$ &   weights $\lambda-Q$ \\
  infinitesimal character $\chi_{\lambda+\rho}$ & infinitesimal character $\chi_{\lambda+\rho}$ \\
  cyclic module & cyclic module \\[-2mm]
  \parbox[c]{4cm}{\begin{center} weight $\mfrak{h}$-module with \\ finite-dimensional \\  weight spaces \end{center}} &  \parbox[c]{4cm}{\begin{center} weight $\mfrak{h}$-module with \\ infinite-dimensional \\  weight spaces \end{center}} \\[-4mm]
  \parbox[c]{4cm}{\begin{center}  $\theta$-Gelfand-Tsetlin module with \\ finite $\Gamma_{\theta}$-multiplicities \\   \end{center}} &  \parbox[c]{4cm}{\begin{center}  $\theta$-Gelfand-Tsetlin module  with \\ finite $\Gamma_{\theta}$-multiplicities  \end{center}}
\end{tabular}
\caption{Comparison of Verma modules and Gelfand--Tsetlin modules}  \label{tab:comparision}
\end{table}

The content of our article goes as follows. In Section \ref{sec:Gelfand-Tsetlin}, we recall a general notion of Harish-Chandra modules and define a class of partial Gelfand--Tsetlin modules for an arbitrary complex simple finite-dimensional Lie algebra $\mfrak{g}$. We also introduce various tensor categories of partial Gelfand--Tsetlin modules. In Section \ref{sec:geometric realization of Gelfand-Tsetlin}, we use the geometric realization of Verma modules $M^\mfrak{g}_\mfrak{b}(\lambda)$ for $\lambda \in \mfrak{h}^*$ described in \cite{Krizka-Somberg2017} to obtain a geometric realization of a new family of $\Gamma_{\theta}$-Gelfand--Tsetlin modules $W^\mfrak{g}_\mfrak{b}(\lambda)$ for $\lambda \in \mfrak{h}^*$ called $\theta$-Gelfand--Tsetlin modules. We also discuss their properties, which are summarized in Table \ref{tab:comparision}. Our main result is the following theorem.
\medskip

\noindent {\bf Theorem.} {\emph{Let $\mfrak{g}$ be a complex simple finite-dimensional Lie algebra with a Cartan subalgebra $\mfrak{h}$ and the corresponding triangular decomposition $\mfrak{g}= \widebar{\mfrak{n}} \oplus \mfrak{h} \oplus \mfrak{n}$. Let $\theta$ be the maximal root in $\Delta^+$ and $\Delta^+_{\theta}=\Delta^+\setminus \{\theta\}$. If $\lambda \in \mfrak{h}^*$, then the vector space $\C[\partial_{x_\alpha},\alpha \in \Delta^+_\theta,\, x_\theta]$ has the structure of a cyclic $\theta$-Gelfand--Tsetlin $\mfrak{g}$-module and every element $a \in \mfrak{g}$ acts as the following differential operator
\begin{align*}
a \mapsto -\sum_{\alpha \in \Delta^+}\bigg[{\ad(u(x))e^{\ad(u(x))} \over e^{\ad(u(x))}-\id}\,(e^{-\ad(u(x))}a)_{\widebar{\mfrak{n}}}\bigg]_\alpha \partial_{x_\alpha} + (\lambda+2\rho)((e^{-\ad(u(x))}a)_\mfrak{h}),
\end{align*}
where $[a]_\alpha$ denotes the $\alpha$-th coordinate
of $a \in \widebar{\mfrak{n}}$ with respect to a root basis $\{f_\alpha;\, \alpha \in \Delta^+\}$
of $\widebar{\mfrak{n}}$, $x=\{x_\alpha;\, \alpha \in \Delta^+\}$  are the corresponding linear coordinate functions on $\widebar{\mfrak{n}}$ and $u(x)=\sum_{\alpha \in \Delta^+} x_\alpha f_\alpha$.
 Further $a_{\widebar{\mfrak{n}}}$ and $a_\mfrak{h}$ are $\widebar{\mfrak{n}}$-part and $\mfrak{h}$-part of $a \in \mfrak{g}$ with respect to the trinagular decomposition $\mfrak{g}=\widebar{\mfrak{n}} \oplus \mfrak{h} \oplus \mfrak{n}$.}}
\medskip

Furthermore, we establish the simplicity criterion for $W^\mfrak{g}_\mfrak{b}(\lambda)$ in the case of  the Lie algebra $\mfrak{sl}(3,\C)$ (Theorem \ref{thm:irreducibility sl(3,C)}). Let us note that the notions of Gelfand--Tsetlin modules and $\theta$-Gelfand--Tsetlin modules coincide for $\mfrak{sl}(3,\C)$. All such simple modules were classified in  \cite{Futorny-Grantcharov-Ramirez2014} (see also \cite{Futorny-Grantcharov-Ramirez2016}).
  We also describe a relation of Gelfand--Tsetlin modules $W^\mfrak{g}_\mfrak{b}(\lambda)$ to twisted Verma modules $M^\mfrak{g}_\mfrak{b}(\lambda)^w$.
   For the reader’s convenience, we remind in Appendix \ref{app:eigenspace decomposition} several important facts concerning the generalized eigenspace decomposition.

Throughout the article, we use the standard notation $\Z$, $\N$ and $\N_0$ for the set of integers, the set of natural numbers and the set of natural numbers together with zero, respectively.


\section{Gelfand--Tsetlin modules}
\label{sec:Gelfand-Tsetlin}

In this section we introduce a general notion of Harish-Chandra modules and then focus on a particular class of partial Gelfand--Tsetlin modules and their analogues for an arbitrary complex simple finite-dimensional Lie algebra.


\subsection{Harish-Chandra modules and Gelfand--Tsetlin modules}

For a commutative $\C$-algebra $\Gamma$ we denote by $\Hom(\Gamma,\C)$ the set of all characters of $\Gamma$, i.e.\ $\C$-algebra homomorphisms from $\Gamma$ to $\C$. Let us note that if $\Gamma$ is finitely generated, then there is a natural identification between the set $\Hom(\Gamma,\C)$ of all characters of $\Gamma$ and the set $\Specm \Gamma$ of all maximal ideals of $\Gamma$, which corresponds to a complex algebraic variety.

Let $M$ be a $\Gamma$-module. For each $\chi \in \Hom(\Gamma,\C)$  we set
\begin{align}
  M_\chi = \{v \in M;\, (\exists k \in \N)\,(\forall a \in \Gamma)\, (a-\chi(a))^kv=0\}
\end{align}
and call it the $\Gamma$-weight space of $M$ with weight $\chi$. When $M_\chi \neq\{0\}$, we say that $\chi$ is a $\Gamma$-weight of $M$ and the elements of $M_\chi$ are called $\Gamma$-weight vectors with weight $\chi$. If a $\Gamma$-module $M$ satisfies
\begin{align}
  M = \bigoplus_{\chi \in \Hom(\Gamma,\C)} M_\chi,
\end{align}
then we call $M$ a $\Gamma$-weight module. The dimension of the vector space $M_\chi$ will be called the \emph{$\Gamma$-multiplicity} of $\chi$ in $M$.
We say that a left $A$-module $M$ is a \emph{Harish-Chandra module with respect to $\Gamma$} if $M$ is a $\Gamma$-weight $A$-module.

\smallskip

Let $A$ be a $\C$-algebra and let $\Gamma$ be a commutative $\C$-subalgebra of $A$. Following \cite{Drozd-Futorny-Ovsienko1994} we denote by $\mcal{H}(A,\Gamma)$ the category of all
$\Gamma$-weight $A$-modules and by $\mcal{H}_{\rm fin}(A,\Gamma)$ the full subcategory of $\mcal{H}(A,\Gamma)$ consisting of those $\Gamma$-weight $A$-modules $M$ satisfying $\dim M_\chi < \infty$ for all $\chi \in \Hom(\Gamma,\C)$. In the case when the $\C$-algebra $A$ is the universal enveloping algebra $U(\mfrak{g})$ of a complex Lie algebra $\mfrak{g}$, we will use the notation $\mcal{H}(\mfrak{g},\Gamma)$ and $\mcal{H}_{\rm fin}(\mfrak{g},\Gamma)$ instead of $\mcal{H}(U(\mfrak{g}),\Gamma)$ and $\mcal{H}_{\rm fin}f(U(\mfrak{g}),\Gamma)$, respectively.

If $\Gamma$ and \smash{$\widetilde{\Gamma}$} are commutative $\C$-subalgebras of $A$ satisfying \smash{$\Gamma \subset \widetilde{\Gamma}$}, then we have the following obvious inclusions
\begin{align*}
\mcal{H}_{\rm fin}(A,\Gamma) \subset \mcal{H}_{\rm fin}(A,\widetilde{\Gamma}) \subset \mcal{H}(A,\widetilde{\Gamma}) \subset \mcal{H}(A,\Gamma),
\end{align*}
which are strict in general.

Finally, let us recall that a commutative $\C$-subalgebra $\Gamma$ of $A$ is a \emph{Harish-Chandra subalgebra} if the $\Gamma$-bimodule $\Gamma a \Gamma$ is finitely generated both as left and as right $\Gamma$-module for all $a \in A$.
\medskip

\example{If $\mathfrak{g}=\mfrak{gl}(n,\C)$ and $\Gamma$  is the Gelfand--Tsetlin subalgebra  of $U(\mathfrak g)$ then objects of $\mcal{H}(\mfrak{g},\Gamma)$ are the Gelfand--Tsetlin modules studied by many authors.}

Let us note that the $\Gamma$-weight spaces in the case of the Gelfand--Tsetlin subalgebra $\Gamma$ of $U(\mathfrak{gl}(n,\C))$ are always finite-dimensional in simple $\mfrak{gl}(n,\C)$-modules. This need not be the case for other commutative $\C$-subalgebras of $U(\mfrak{gl}(n,\C))$. However, in all known to us cases the following is true. We state it as a conjecture in general.
\medskip

\conjecture{Let $\mathfrak g$ be a complex simple finite-dimensional Lie algebra. If $M \in \mcal{H}(\mathfrak{g}, \Gamma)$ is a simple $\mfrak{g}$-module and $\dim M_\chi <\infty$ for some $\Gamma$-weight $\chi$ of $M$, then $M \in \mcal{H}_{\rm fin}(\mathfrak{g},\Gamma)$.}

\remark{If $M$ is a simple $\mathfrak{g}$-module such that  $M_{\chi}\neq \{0\}$ for some $\chi \in \Hom(\Gamma,\C)$, then it is not clear (and probably not true in general) whether $M \in \mcal{H}(\mathfrak{g},\Gamma)$.}
\vspace{-2mm}


\subsection{$\theta$-Gelfand--Tsetlin modules}

In this section we introduce a class of partial Gelfand--Tsetlin modules for an arbitrary complex
simple finite-dimensional Lie algebra.
\medskip

Let us consider a complex simple finite-dimensional Lie algebra $\mfrak{g}$ and let $\mfrak{h}$ be a Cartan subalgebra of $\mfrak{g}$. We denote by $\Delta$ the root system of $\mfrak{g}$ with respect to $\mfrak{h}$, by $\Delta^+$ a positive root system in $\Delta$ and by $\Pi \subset \Delta$ the set of simple roots. Further, we denote by $\theta$ the maximal root of $\mfrak{g}$ (by definition, the highest weight of its adjoint representation) and by $h_\alpha \in \mfrak{h}$ the coroot corresponding to a root $\alpha \in \Delta$.
\medskip

\definition{Let  $\mfrak{s}_\theta$ be the Lie subalgebra of $\mfrak{g}$ generated by an $\mfrak{sl}(2,\C)$-triple $(e_\theta,h_\theta,f_\theta)$ associated to the maximal root $\theta$. Then we define the commutative $\C$-subalgebra $\Gamma_\theta$ of $U(\mfrak{g})$ as a $\C$-subalgebra of $U(\mfrak{g})$ generated by the Cartan subalgebra $\mfrak{h}$ and by the center $\mfrak{Z}(\mfrak{s}_\theta)$ of $U(\mfrak{s}_\theta)$.}

The center $\mfrak{Z}(\mfrak{s}_\theta)$ of $U(\mfrak{s}_\theta)$ is freely generated by the quadratic Casimir element $\Cas_\theta$ defined by
\begin{align}
   \Cas_\theta = e_\theta f_\theta + f_\theta e_\theta + {\textstyle {1 \over 2}} h_\theta^2.
\end{align}
Hence, the $\C$-algebra $\Gamma_\theta$ is freely generated by the coroots $h_\alpha$, $\alpha \in \Pi$, and by the Casimir element $\Cas_\theta$. As a consequence, we obtain that the complex algebraic variety corresponding to $\Hom(\Gamma_\theta,\C)$ and $\Specm \Gamma_\theta$ is isomorphic to $\C^{\rank(\mfrak{g})+1}$. In the next, we will focus on the category $\mcal{H}(\mathfrak{g}, \Gamma_\theta)$. The objects of this category will be called \emph{$\theta$-Gelfand--Tsetlin modules}.
\medskip

\proposition{
\begin{enum}
\item[(i)] The commutative $\C$-subalgebra $\Gamma_\theta$ of $U(\mfrak{g})$ is a Harish-Chandra subalgebra of $U(\mfrak{g})$.
\item[(ii)] Let $M$ be a simple $\mfrak{g}$-module such that $M_\chi \neq \{0\}$ for some $\chi \in \Hom(\Gamma_\theta,\C)$. Then we have $M\in \mcal{H}(\mathfrak{g},\Gamma_\theta)$.
\end{enum}
}

\proof{To prove the first statement it is sufficient to check it for the generators of the Lie algebra $\mfrak{g}$. Let $a \in \mfrak{g}$ be such a generator. We may assume that $a\in \mfrak{g}_{\alpha}$ for some $\alpha \in \Delta$. If $\alpha = \pm \theta$, then we have $[\Cas_\theta,a]=0$. For $\alpha \neq \pm \theta$, let us consider the $\mfrak{s}_\theta$-module
\begin{align*}
\mfrak{g}_{\alpha,\theta}=\bigoplus_{j \in \Z} \mfrak{g}_{\alpha + j\theta}.
\end{align*}
Since $\mfrak{g}_{\alpha,\theta}$ is a simple finite-dimensional $\mfrak{s}_\theta$-module, we obtain that $[\Cas_\theta,a]=c_\alpha a$, where $c_\alpha \in \C$, which together with $[h,a]=\alpha(h) a$ implies the first and also the second statement. }

We define the full subcategories $\mcal{I}_{\theta, e}$ and $\mcal{I}_{\theta, f}$ of the category of weight (with respect to $\mfrak{h}$) $\mfrak{g}$-modules consisting of those $\mfrak{g}$-modules on which $e_\theta$ and $f_\theta$ is locally nilpotent, respectively. Also we set $\mcal{I}_\theta=\mcal{I}_{\theta,e} \cap \mcal{I}_{\theta,f}$.
\medskip

\proposition{
 \begin{enum}
\item[(i)] The categories $\mcal{I}_{\theta,e}$, $\mcal{I}_{\theta,f}$ and $\mcal{I}_\theta$ are full subcategories of $\mcal{H}(\mathfrak{g}, \Gamma_\theta)$.
\item[(ii)] Let $M$ be a simple weight $\mfrak{g}$-module such that there exists an integer $n \in \N$ and a nonzero vector $v \in M$ satisfying $e_\theta^n v=0$ (respectively, $f_\theta^n v=0$). Then we have
   $M \in \mcal{I}_{\theta,e}$ (respectively, $M \in  \mcal{I}_{\theta,f}$), and hence  $M\in \mcal{H}(\mathfrak{g}, \Gamma_{\theta})$.
 \end{enum}
}

\proof{If $M \in \mcal{I}_{\theta,e}$, then $\Cas_\theta$ is locally finite on $M$. This implies the first statement. The second statement follows from the simplicity of $M$ and item (i).}

We will denote $\mcal{I}_{\theta,e,{\rm fin}}= \mcal{I}_{\theta, e}\cap \mcal{H}_{\rm fin}(\mathfrak{g}, \Gamma_\theta)$, and similarly $\mcal{I}_{\theta, f,{\rm fin}}$, $\mcal{I}_{\theta,{\rm fin}}$. Let $\Gamma$ be a commutative $\C$-subalgebra of $U(\mfrak{g})$ containing $\mfrak{h}$ and $\mfrak{Z}(\mfrak{s}_\theta)$. Since  $\mcal{H}_{\rm fin}(\mathfrak{g},\Gamma_\theta) \subset \mcal{H}_{\rm fin}(\mathfrak{g}, \Gamma)$, we immediately have the following statement.
\medskip

\corollary{The categories $\mcal{I}_{\theta, e, {\rm fin}}$, $\mcal{I}_{\theta, f,{\rm fin}}$ and $\mcal{I}_{\theta,{\rm fin}}$ are full subcategories of $\mcal{H}_{\rm fin}(\mathfrak{g},\Gamma)$.}

One wonders how large all these categories are. Clearly, $\mcal{I}_{\theta,{\rm fin}}$ contains all finite-dimensional $\mathfrak{g}$-modules; $\mcal{I}_{\theta, e,{\rm fin}}$ contains the whole
category $\mcal{O}$ and all parabolically induced $\mfrak{g}$-modules when the Levi factor $\mathfrak{l}$ of the parabolic subalgebra $\mathfrak{p}=\mathfrak{l} \oplus \mathfrak{u}$ of $\mathfrak{g}$ does not contain  $e_\theta$, $f_\theta$ and $e_\theta \in \mathfrak{u}$.
Moreover, it was shown in \cite{Futorny-Grantcharov-Ramirez2014} that certain simple non-parabolically induced Gelfand--Tsetlin $\mfrak{sl}(3,\C)$-modules with respect to the Gelfand--Tsetlin subalgebra belong to the categories $\mcal{I}_{\theta, e,{\rm fin}}$, $\mcal{I}_{\theta, f,{\rm fin}}$ and $\mcal{I}_{\theta,{\rm fin}}$. Such modules have infinite-dimensional weight spaces. Furthermore, using the parabolic induction from such
simple Gelfand--Tsetlin $\mfrak{sl}(3,\C)$-modules  one can construct  simple $\mathfrak g$-modules respectively in categories $\mcal{I}_{\theta, e}$, $\mcal{I}_{\theta, f}$ and $\mcal{I}_\theta$ for any $\mathfrak{g}$ of rank $\rank(\mfrak{g}) \geq 3$ (in order to have $\mfrak{sl}(3,\C)+\mathfrak{h}$ as a Levi factor of some parabolic subalgebra).

Let us note that $\mcal{H}(\mathfrak{g},\Gamma_\theta)$, $\mcal{H}_{\rm fin}(\mathfrak{g}, \Gamma_\theta)$, $\mcal{I}_{\theta, e,{\rm fin}}$, $\mcal{I}_{\theta, f,{\rm fin}}$ and $\mcal{I}_{\theta,{\rm fin}}$ are not tensor categories, except the case $\mathfrak{g}=\mfrak{sl}(2,\C)$ when $\mcal{I}_{\theta, e, {\rm fin}}$ (respectively, $\mcal{I}_{\theta, f,{\rm fin}}$) consists of extensions of highest (respectively, lowest) weight $\mfrak{g}$-modules and their directs sums and $\mcal{I}_{\theta,{\rm fin}}$ consists of sums of finite-dimensional $\mfrak{g}$-modules. Nevertheless, we have the following statement as a consequence of the definitions.

\medskip

\theorem{\emph{The categories $\mcal{I}_{\theta, e}$, $\mcal{I}_{\theta, f}$ and $\mcal{I}_\theta$ are tensor categories.}}

\vspace{-2mm}

\section{Geometric construction of $\theta$-Gelfand--Tsetlin modules}\label{sec:geometric realization of Gelfand-Tsetlin}


\subsection{Geometric realization of Verma modules}
Let us consider a complex simple Lie algebra $\mfrak{g}$ and let $\mfrak{h}$ be a Cartan subalgebra of $\mfrak{g}$. We denote by $\Delta$ the root system of $\mfrak{g}$ with respect to $\mfrak{h}$, by $\Delta^+$ a positive root system in $\Delta$ and by $\Pi \subset \Delta$ the set of simple roots. The standard Borel subalgebra $\mfrak{b}$ of $\mfrak{g}$ is defined through $\mfrak{b} = \mfrak{h} \oplus \mfrak{n}$ with the nilradical $\mfrak{n}$ and the opposite nilradical $\widebar{\mfrak{n}}$ given by
\begin{align}
  \mfrak{n} = \bigoplus_{\alpha \in \Delta^+} \mfrak{g}_\alpha \qquad \text{and} \qquad \widebar{\mfrak{n}} = \bigoplus_{\alpha \in \Delta^+} \mfrak{g}_{-\alpha}. \label{eq:niradical borel}
\end{align}
Moreover, we have a triangular decomposition
\begin{align}
  \mfrak{g}= \widebar{\mfrak{n}} \oplus \mfrak{h} \oplus \mfrak{n} \label{eq:triangular decomposition}
\end{align}
of the Lie algebra $\mfrak{g}$. Furthermore, we define the height $\htt(\alpha)$ of $\alpha \in \Delta$ by
\begin{align}
  \htt({\textstyle \sum_{i=1}^r} a_i\alpha_i)={\textstyle \sum^r_{i=1}}  a_i,
\end{align}
where $r= \rank(\mfrak{g})$ and $\Pi=\{\alpha_1,\alpha_2,\dots,\alpha_r\}$. If we denote $k=\htt(\theta)$, where $\theta$ is the maximal root of $\mfrak{g}$ (by definition, the highest weight of its adjoint representation), then $\mfrak{g}$ is a $|k|$-graded Lie algebra with respect to the grading given by $\mfrak{g}_i= \bigoplus_{\alpha \in \Delta,\, \htt(\alpha)=i} \mfrak{g}_\alpha$ for $0 \neq i \in \Z$ and $\mfrak{g}_0 = \mfrak{h}$. Moreover, we have
\begin{align}
\widebar{\mfrak{n}}= \mfrak{g}_{-k} \oplus \dots \oplus \mfrak{g}_{-1}, \qquad \mfrak{h} =\mfrak{g}_0, \qquad  \mfrak{n}= \mfrak{g}_1 \oplus \dots \oplus \mfrak{g}_k
\end{align}
together with $\mfrak{g}_{-k}=\mfrak{g}_{-\theta}$ and $\mfrak{g}_k = \mfrak{g}_\theta$. We also set
\begin{align}
  Q = \sum_{\alpha \in \Delta^+} \Z \alpha = \bigoplus_{i=1}^r \Z \alpha_i \qquad \text{and} \qquad Q_+= \sum_{\alpha \in \Delta^+} \N_0 \alpha = \bigoplus_{i=1}^r \N_0 \alpha_i
\end{align}
and call $Q$ the root lattice. Let us note that for any positive root $\alpha \in \Delta^+$, we have $\theta -\alpha \in Q_+$. Furthermore, we introduce the notation $\Delta^+_\theta$ to stand for $\Delta^+ \setminus \{\theta\}$.
\medskip

For $\lambda \in \mfrak{h}^*$ we denote by $\C_\lambda$ the $1$-dimensional representation of $\mfrak{b}$ defined by
\begin{align}
  av=\lambda(a)v
\end{align}
for $a \in \mfrak{h}$, $v \in \C_\lambda \simeq \C$ (as vector spaces) and trivially extended to a representation of $\mfrak{b} = \mfrak{h} \oplus \mfrak{n}$.
\medskip

\definition{Let $\lambda \in \mfrak{h}^*$. Then the Verma module $M^\mfrak{g}_\mfrak{b}(\lambda)$ is the induced module
\begin{align}
  M^\mfrak{g}_\mfrak{b}(\lambda) = \Ind^\mfrak{g}_\mfrak{b} \C_\lambda \equiv U(\mfrak{g}) \otimes_{U(\mfrak{b})} \C_\lambda \simeq U(\widebar{\mfrak{n}}) \otimes_\C \C_\lambda,
\end{align}
where the last isomorphism of $U(\widebar{\mfrak{n}})$-modules follows from  Poincaré--Birkhoff--Witt theorem.}

In the next, we shall need the action of the center $\mfrak{Z}(\mfrak{g})$ of the universal enveloping algebra $U(\mfrak{g})$ on Verma modules $M^\mfrak{g}_\mfrak{b}(\lambda)$ for $\lambda \in \mfrak{h}^*$. Using \eqref{eq:triangular decomposition} together with Poincaré--Birkhof--Witt theorem we obtain
\begin{align}
  U(\mfrak{g}) = U(\mfrak{h}) \oplus (\widebar{\mfrak{n}}\,U(\mfrak{g}) + U(\mfrak{g})\,\mfrak{n}). \label{eq:universal enveloping decomposition}
\end{align}
Let us consider the projection $p\colon U(\mfrak{g}) \rarr U(\mfrak{h})$ with respect to the direct sum decomposition \eqref{eq:universal enveloping decomposition} and the $\C$-algebra automorphism $f \colon U(\mfrak{h}) \rarr U(\mfrak{h})$ defined by
\begin{align}
  f(h) = h-\rho(h)1
\end{align}
for all $h \in \mfrak{h}$, where $\rho \in \mfrak{h}^*$ is the Weyl vector given through
\begin{align}
 \rho =  {1 \over 2} \sum_{\alpha \in \Delta^+} \alpha. \label{eq:weyl vector}
\end{align}
We denote by $\gamma \colon \mfrak{Z}(\mfrak{g}) \rarr U(\mfrak{h})$ the restriction of $f \circ p$ to the center $\mfrak{Z}(\mfrak{g})$ of $U(\mfrak{g})$.
The mapping $\gamma \colon \mfrak{Z}(\mfrak{g}) \rarr U(\mfrak{h})$ is a homomorphism of associative $\C$-algebras and does not depend on the choice of the positive root system $\Delta^+$. It is called the Harish-Chandra homomorphism.

A $\C$-algebra homomorphism from $\mfrak{Z}(\mfrak{g})$ to $\C$ is called a central character. For each $\lambda \in \mfrak{h}^*$ we define a central character $\chi_\lambda \colon \mfrak{Z}(\mfrak{g}) \rarr \C$ by
\begin{align}
  \chi_\lambda(z)=(\gamma(z))(\lambda)
\end{align}
for $z \in \mfrak{Z}(\mfrak{g})$, where we identify $U(\mfrak{h}) \simeq S(\mfrak{h})$ with the $\C$-algebra of polynomial functions on $\mfrak{h}^*$. The action of the center $\mfrak{Z}(\mfrak{g})$ of $U(\mfrak{g})$ on Verma modules $M^\mfrak{g}_\mfrak{b}(\lambda)$ for $\lambda \in \mfrak{h}^*$ is then given by
\begin{align}
  zv=\chi_{\lambda+\rho}(z)v \label{eq:center action Verma}
\end{align}
for all $z \in \mfrak{Z}(\mfrak{g})$ and $v \in M^\mfrak{g}_\mfrak{b}(\lambda)$.
\medskip

We briefly describe the geometric realization of Verma modules given in \cite{Krizka-Somberg2017}. This nice geometric realization enables us to construct new interesting simple modules for the Lie algebra $\mfrak{g}$.

Let $\{f_\alpha;\, \alpha \in \Delta^+\}$ be a basis of the opposite nilradical $\widebar{\mfrak{n}}$. We denote by $x=\{x_\alpha;\, \alpha \in \Delta^+\}$ the corresponding linear coordinate functions on $\widebar{\mfrak{n}}$ with respect to the basis $\{f_\alpha;\, \alpha \in \Delta^+\}$ of $\widebar{\mfrak{n}}$. Then we have a canonical isomorphism $\C[\widebar{\mfrak{n}}] \simeq \C[x]$ of $\C$-algebras, 
and the Weyl algebra $\eus{A}^\mfrak{g}_{\widebar{\mfrak{n}}}$ of the vector space $\widebar{\mfrak{n}}$ is generated by $\{x_\alpha, \partial_{x_\alpha};\, \alpha \in \Delta^+\}$ together with the canonical commutation relations. For $\lambda \in \mfrak{h}^*$ there is a homomorphism
\begin{align}
  \pi_\lambda \colon \mfrak{g} \rarr \eus{A}^\mfrak{g}_{\widebar{\mfrak{n}}}
\end{align}
of Lie algebras given through
\begin{align}
\pi_\lambda(a)= -\sum_{\alpha \in \Delta^+}\bigg[{\ad(u(x))e^{\ad(u(x))} \over e^{\ad(u(x))}-\id}\,(e^{-\ad(u(x))}a)_{\widebar{\mfrak{n}}}\bigg]_\alpha \partial_{x_\alpha} + (\lambda+\rho)((e^{-\ad(u(x))}a)_\mfrak{h}) \label{eq:pi action general}
\end{align}
for all $a \in \mfrak{g}$, where $[a]_\alpha$ denotes the $\alpha$-th coordinate
of $a \in \widebar{\mfrak{n}}$ with respect to the basis $\{f_\alpha;\, \alpha \in \Delta^+\}$
of $\widebar{\mfrak{n}}$, further $a_{\widebar{\mfrak{n}}}$ and $a_\mfrak{h}$ are $\widebar{\mfrak{n}}$-part and $\mfrak{h}$-part of $a \in \mfrak{g}$ with respect to the trinagular decomposition $\mfrak{g}=\widebar{\mfrak{n}} \oplus \mfrak{h} \oplus \mfrak{n}$, and finally the element $u(x) \in \C[\widebar{\mfrak{n}}] \otimes_\C \mfrak{g}$ is given by
\begin{align}
u(x)=\sum_{\alpha \in \Delta^+} x_\alpha f_\alpha.
\end{align}
Let us note that $\C[\widebar{\mfrak{n}}] \otimes_\C \mfrak{g}$ has the natural structure of a Lie algebra. Hence, we have a well-defined $\C$-linear mapping $\ad(u(x)) \colon \C[\widebar{\mfrak{n}}] \otimes_\C \mfrak{g} \rarr \C[\widebar{\mfrak{n}}] \otimes_\C \mfrak{g}$.

In particular, we have
\begin{align}
  \pi_\lambda(a)= - \sum_{\alpha \in\Delta^+} \bigg[{\ad(u(x)) \over e^{\ad(u(x))} - \id}\,a\bigg]_\alpha \partial_{x_\alpha} \label{eq:pi action nilradical}
\end{align}
for $a \in \widebar{\mfrak{n}}$ and
\begin{align}
  \pi_\lambda(a)= \sum_{\alpha \in \Delta^+} [\ad(u(x))a]_\alpha \partial_{x_\alpha} + (\lambda+\rho)(a) \label{eq:pi action cartan}
\end{align}
for $a \in \mfrak{h}$. The Verma module $M^\mfrak{g}_\mfrak{b}(\lambda)$ for $\lambda \in \mfrak{h}^*$ is then realized as $\eus{A}^\mfrak{g}_{\widebar{\mfrak{n}}}/\eus{I}_{\rm V}$, where
$\eus{I}_{\rm V}$ is the left ideal of $\eus{A}^\mfrak{g}_{\widebar{\mfrak{n}}}$ defined by $\eus{I}_{\rm V}=(x_\alpha,\, \alpha \in \Delta^+)$, see e.g.\ \cite{Krizka-Somberg2017}. The $\mfrak{g}$-module structure on $\eus{A}^\mfrak{g}_{\widebar{\mfrak{n}}}/\eus{I}_{\rm V}$ is given through the homomorphism
$\pi_{\lambda+\rho} \colon \mfrak{g} \rarr \eus{A}^\mfrak{g}_{\widebar{\mfrak{n}}}$ of Lie algebras.

In the next, we change the left ideal $\eus{I}_{\rm V}$ of the Weyl algebra $\eus{A}^\mfrak{g}_{\widebar{\mfrak{n}}}$ to produce a class of $\theta$-Gelfand--Tsetlin modules. These $\mfrak{g}$-modules will be the subject of our interest.


\subsection{Geometric realization of $\theta$-Gelfand--Tsetlin modules}

We need to introduce a suitable basis of the Lie algebra $\mfrak{g}$. Let $\{e_\alpha;\, \alpha \in \Delta^+\}$ be a root basis of the nilradical $\mfrak{n}$ and let $\{f_\alpha;\, \alpha \in \Delta^+\}$ be a root basis of the opposite nilradical $\widebar{\mfrak{n}}$. If we denote by $h_\alpha \in \mfrak{h}$ the coroot corresponding to a root $\alpha \in \Delta^+$, then $\{h_\alpha;\, \alpha\in \Pi\}$ is a basis of the Cartan subalgebra $\mfrak{h}$.

Let us note that the Lie subalgebra of $\mfrak{g}$ generated by $\{e_\alpha,h_\alpha,f_\alpha\}$ for $\alpha \in \Delta^+$ is isomorphic to $\mfrak{sl}(2,\C)$.
\medskip

\definition{Let $\lambda \in \mfrak{h}^*$. Then the $\mfrak{g}$-module $W^\mfrak{g}_\mfrak{b}(\lambda)$ is defined by
\begin{align}
W^\mfrak{g}_\mfrak{b}(\lambda) \simeq \eus{A}^\mfrak{g}_{\widebar{\mfrak{n}}}/ \eus{I}_{\rm GT} \simeq
\C[\partial_{x_\alpha},\alpha \in \Delta^+_\theta,\, x_\theta], \quad \eus{I}_{\rm GT}=(x_\alpha, \alpha \in \Delta^+_\theta,\partial_{x_\theta}),
\end{align}
where the $\mfrak{g}$-module structure on $\eus{A}^\mfrak{g}_{\widebar{\mfrak{n}}}/ \eus{I}_{\rm GT}$ is given through the homomorphism $\pi_{\lambda+\rho} \colon  \mfrak{g}\rarr \eus{A}^\mfrak{g}_{\widebar{\mfrak{n}}}$ of Lie algebras.}

\remark{Although, the definition of $W^\mfrak{g}_\mfrak{b}(\lambda)$ for $\lambda \in \mfrak{h}^*$ makes sense for any complex simple finite-dimensional Lie algebra $\mfrak{g}$, we exclude the case $\mfrak{g}=\mfrak{sl}(2,\C)$. For $\mfrak{g}=\mfrak{sl}(2,\C)$ the $\mfrak{g}$-module $W^\mfrak{g}_\mfrak{b}(\lambda)$ is isomorphic to the contragredient Verma module \smash{$M^\mfrak{g}_{\widebar{\mfrak{b}}}(\lambda+2\rho)^*$} for the opposite standard Borel subalgebra $\widebar{\mfrak{b}} = \mfrak{h} \oplus \widebar{\mfrak{n}}$.}

Let $M$ be an $\mfrak{h}$-module. For each $\mu \in \mfrak{h}^*$ we set
\begin{align}
  M_\mu = \{v \in M;\, (\forall h \in \mfrak{h})\, hv=\mu(h)v\}
\end{align}
and call it the weight space of $M$ with weight $\mu$. When $M_\mu \neq \{0\}$, we say that $\mu$ is a weight of $M$ and the elements of $M_\mu$ are called weight vectors with weight $\mu$. If an $\mfrak{h}$-module $M$ satisfies
\begin{align}
  M = \bigoplus_{\mu\in \mfrak{h}^*} M_\mu,
\end{align}
then we call $M$ a weight module.
\medskip

Therefore, the first natural question is whether $W^\mfrak{g}_\mfrak{b}(\lambda)$ is a weight module and how the action of the center $\mfrak{Z}(\mfrak{g})$ of $U(\mfrak{g})$ on $W^\mfrak{g}_\mfrak{b}(\lambda)$ looks.
\medskip

\lemma{Let $\lambda \in \mfrak{h}^*$. Then we have
\begin{align} \label{eq:action weight}
\begin{aligned}
  \pi_\lambda(h)&= \sum_{\alpha \in \Delta^+} \alpha(h) x_\alpha \partial_{x_\alpha} + (\lambda+\rho)(h) \\
  &=\sum_{\alpha \in \Delta^+_\theta} \alpha(h) \partial_{x_\alpha} x_\alpha + \theta(h) x_\theta \partial_{x_\theta} + (\lambda-\rho+\theta)(h)
\end{aligned}
\end{align}
for all $h \in \mfrak{h}$.}

\proof{For $h \in \mfrak{h}$ we have
\begin{align*}
  [u(x),h]= \sum_{\alpha \in \Delta^+} x_\alpha [f_\alpha,h] = \sum_{\alpha \in \Delta^+} \alpha(h) x_\alpha f_\alpha,
\end{align*}
which together with \eqref{eq:pi action cartan} enables us to write
\begin{align*}
  \pi_\lambda(h)&= \sum_{\alpha \in \Delta^+} \alpha(h) x_\alpha \partial_{x_\alpha} + (\lambda+\rho)(h) \\
  &=\sum_{\alpha \in \Delta^+_\theta} \alpha(h) \partial_{x_\alpha} x_\alpha + \theta(h) x_\theta \partial_{x_\theta} + (\lambda-\rho+\theta)(h)
\end{align*}
for all $h \in \mfrak{h}$, where we used \eqref{eq:weyl vector} in the second equality.}

We denote by $\Z^{\Delta^+}$\! and $\Z^\Pi$ the set of all functions from $\Delta^+$ to $\Z$ and from $\Pi$ to $\Z$, respectively. Since we have $\Pi \subset \Delta^+$, an element of $\Z^\Pi$ will be also regarded as an element of $\Z^{\Delta^+}$\! extended by $0$ on $\Delta^+ \setminus \Pi$. A similar notation is introduced for $\N_0^{\Delta^+}$\! and $\N_0^\Pi$.

For \smash{$a \in \N_0^{\Delta^+}$}\! we define by
\begin{align}
   w_{\lambda,a} =  \prod_{\alpha \in \Delta^+_\theta} \partial_{x_\alpha}^{a_\alpha} x_\theta^{a_\theta}
\end{align}
a vector in $W^\mfrak{g}_\mfrak{b}(\lambda)$. The subset $\{w_{\lambda,a};\, a \in \N_0^{\Delta^+}\!\}$ of $W^\mfrak{g}_\mfrak{b}(\lambda)$ forms a basis of $W^\mfrak{g}_\mfrak{b}(\lambda)$. As $w_{\lambda,a}$ is a weight vector with weight
\begin{align}
   \mu_{\lambda,a}  = -  \sum\limits_{\alpha \in \Delta^+_\theta} a_\alpha \alpha + a_\theta \theta  + \lambda + \theta \label{eq:weights}
\end{align}
for all $a \in \N_0^{\Delta^+}$\!, as easily follows from \eqref{eq:action weight}, we obtain that $W^\mfrak{g}_\mfrak{b}(\lambda)$ is a weight module. However, we need a basis of weight spaces. Since any positive root $\alpha \in \Delta^+$ can be expressed as
\begin{align}
  \alpha = \sum_{i=1}^r m_{\alpha,i} \alpha_i, \label{eq:simple root decomposition}
\end{align}
where $m_{\alpha,i} \in \N_0$ for all $i=1,2,\dots,r$ and $\Pi=\{\alpha_1, \alpha_2, \dots, \alpha_r\}$, we define $t_\alpha \in \Z^{\Delta^+}$\! for $\alpha \in \Delta^+ \setminus \Pi$ by
\begin{align}
  t_\theta(\beta) = \begin{cases}
    m_{\theta,i} &  \text{for $\beta=\alpha_i$}, \\
    1 & \text{for $\beta=\theta$}, \\
    0 & \text{for $\beta \neq \theta$ and $\beta \notin \Pi$}
  \end{cases}
\end{align}
and
\begin{align}
  t_\alpha(\beta) = \begin{cases}
    -m_{\alpha,i} &  \text{for $\beta=\alpha_i$}, \\
    1 & \text{for $\beta=\alpha$}, \\
    0 & \text{for $\beta \neq \alpha$ and $\beta \notin \Pi$}
  \end{cases}
\end{align}
for $\alpha \neq \theta$. Furthermore, we define subsets $\Lambda_+$ and $\Lambda_+^\theta$ of $\Z^{\Delta^+}$\! by
\begin{align}
  \Lambda_+ = \{{\textstyle \sum_{\alpha \in \Delta^+ \setminus \Pi} n_\alpha t_\alpha};\, n_\alpha \in \N_0\ \text{for all $\alpha \in \Delta^+ \setminus \Pi$}\}
\end{align}
and
\begin{align}
  \Lambda_+^\theta = \{{\textstyle \sum_{\alpha \in \Delta^+_\theta \setminus \Pi} n_\alpha t_\alpha};\, n_\alpha \in \N_0\ \text{for all $\alpha \in \Delta^+_\theta \setminus \Pi$}\},
\end{align}
respectively.

\medskip

\lemma{\label{lem:mu weights}
Let $\lambda \in \mfrak{h}^*$ and $a,b \in \Z^{\Delta^+}$\!. Then $\mu_{\lambda,a} = \mu_{\lambda,b}$ if and only if
\begin{align}
  b = a + \sum_{\alpha \in \Delta^+ \setminus \Pi} n_\alpha t_\alpha,
\end{align}
where $n_\alpha \in \Z$ for all $\alpha \in \Delta^+ \setminus \Pi$.}

\proof{Let $\lambda \in \mfrak{h}^*$ and $a,b \in \Z^{\Delta^+}$\!. If \smash{$b=a+ \sum_{\alpha \in \Delta^+ \setminus \Pi} n_\alpha t_\alpha$}, where $n_\alpha \in \Z$ for $\alpha \in \Delta^+ \setminus \Pi$, then we have $\mu_{\lambda,a}=\mu_{\lambda,b}$. On the other hand, let us assume that $\mu_{\lambda,a}=\mu_{\lambda,b}$. Then we set $n_\alpha = b_\alpha-a_\alpha$ for $\alpha \in \Delta^+ \setminus \Pi$ and define \smash{$c=a + \sum_{\alpha \in \Delta^+ \setminus \Pi} n_\alpha t_\alpha$}. Hence, we have $c_\alpha = a_\alpha + n_\alpha = b_\alpha$ for $\alpha \in \Delta^+ \setminus \Pi$ and \smash{$c_{\alpha_i} = a_{\alpha_i} - \sum_{\alpha \in \Delta^+_\theta \setminus \Pi} n_\alpha m_{\alpha,i} + n_\theta m_{\theta,i}$} for $i=1,2,\dots,r$. Further, since $\mu_{\lambda,a}=\mu_{\lambda,b}$, we may write
\begin{align*}
  \mu_{\lambda,b}-\mu_{\lambda,a} &= \sum_{\alpha \in \Delta^+_\theta} (a_\alpha-b_\alpha) \alpha - (a_\theta-b_\theta)\theta \\
  &= \sum_{i=1}^r (a_{\alpha_i}-b_{\alpha_i})\alpha_i + \sum_{i=1}^r \sum_{\alpha \in \Delta^+_\theta \setminus \Pi} (a_\alpha - b_\alpha)m_{\alpha,i}\alpha_i - \sum_{i=1}^r (a_\theta - b_\theta)m_{\theta,i}\alpha_i \\
  &= \sum_{i=1}^r \!\Big(a_{\alpha_i}-b_{\alpha_i} + \sum_{\alpha \in \Delta^+_\theta \setminus \Pi} (a_\alpha - b_\alpha)m_{\alpha,i} - (a_\theta - b_\theta)m_{\theta,i}\Big)\alpha_i \\
  &= \sum_{i=1}^r \!\Big(a_{\alpha_i}-b_{\alpha_i} - \sum_{\alpha \in \Delta^+_\theta \setminus \Pi} n_\alpha m_{\alpha,i} + n_\theta m_{\theta,i}\Big)\alpha_i=0,
\end{align*}
where we used \eqref{eq:simple root decomposition}. As the set $\{\alpha_1,\alpha_2,\dots,\alpha_r\}$ forms a basis of $\mfrak{h}^*$, we get
\begin{align*}
 a_{\alpha_i}-b_{\alpha_i} - \sum_{\alpha \in \Delta^+_\theta \setminus \Pi} n_\alpha m_{\alpha,i} + n_\theta m_{\theta,i}=0,
\end{align*}
which implies that $c_{\alpha_i}=b_{\alpha_i}$ for $i=1,2,\dots,r$. Hence, we have $c_\alpha= b_\alpha$ for all $\alpha \in \Delta^+$ and we are done.}

\proposition{\label{prop:weights description}
Let $\lambda  \in \mfrak{h}^*$. Then the set of all weights of $W^\mfrak{g}_\mfrak{b}(\lambda)$ is $\{\mu_{\lambda,a};\, a \in \Z^\Pi\}$. Moreover, we have that $W^\mfrak{g}_\mfrak{b}(\lambda)$ is a weight $\mfrak{h}$-module with infinite dimensional weight spaces and the subset
\begin{align}
  \{w_{\lambda,a+n_a t_\theta+t};\, t \in \Lambda_+,\, a+n_a t_\theta +t \in \N_0^{\Delta^+}\!\}
\end{align}
of $W^\mfrak{g}_\mfrak{b}(\lambda)_{\mu_{\lambda,a}}$ forms a basis of $W^\mfrak{g}_\mfrak{b}(\lambda)_{\mu_{\lambda,a}}$ for all $a \in \Z^\Pi$, where $n_a \in \N_0$ is the smallest non-negative integer satisfying \smash{$a+n_at_\theta \in \N_0^{\Delta^+}$}\!.}

\proof{From the previous considerations we know that $W^\mfrak{g}_\mfrak{b}(\lambda)$ is a weight $\mfrak{h}$-module and that all weighs are of the form $\mu_{\lambda,a}$ for some $a \in \Z^\Pi$ as follows from \eqref{eq:weights}. Now, let us consider $a \in \Z^\Pi$. Since $\theta - \alpha \in Q_+$ for all $\alpha \in \Delta^+$, i.e.\ $m_{\theta,i} \neq 0$ for all $i=1,2,\dots,r$, there exists the smallest non-negative integer $n_a \in \N_0$ such that \smash{$a+n_a t_\theta \in \N_0^{\Delta^+}$}\!. Moreover, we have $\mu_{\lambda,a} = \mu_{\lambda,a+n_at_\theta}$, which implies that $\mu_{\lambda,a}$ is a weight of $W^\mfrak{g}_\mfrak{b}(\lambda)$, since $w_{\lambda,a+n_at_\theta}$ is a weight vector with weight $\mu_{\lambda,a}$. Further, since \smash{$a+(n_a+n)t_\theta \in \N_0^{\Delta^+}$}\! and $\mu_{\lambda,a} = \mu_{\lambda,a+(n_a+n)t_\theta}$ for all $n \in \N_0$, we obtain that $W^\mfrak{g}_\mfrak{b}(\lambda)_{\mu_{\lambda,a}}$ contains the subset $\{w_{\lambda,a+(n_a+n)t_\theta};\, n \in \N_0\}$, which implies that the weight space $W^\mfrak{g}_\mfrak{b}(\lambda)_{\mu_{\lambda,a}}$ is  infinite dimensional.

Let us assume that $w_{\lambda,b} \in W^\mfrak{g}_\mfrak{b}(\lambda)_{\mu_{\lambda,a}}$ for $a \in \Z^\Pi$ and \smash{$b \in \N_0^{\Delta^+}$}\!. Then we have $\mu_{\lambda,a}=\mu_{\lambda,b}$. By Lemma \ref{lem:mu weights} we get
\begin{align*}
  b = a + \sum_{\alpha\in \Delta^+ \setminus \Pi} n_\alpha t_\alpha,
\end{align*}
where $n_\alpha \in \Z$ for all $\alpha \in \Delta^+ \setminus \Pi$. Therefore, we have $b_\alpha = a_\alpha +n_\alpha = n_\alpha$ for $\alpha \in \Delta^+ \setminus \Pi$ and \smash{$b_{\alpha_i}= a_{\alpha_i}  +n_\theta m_{\theta,i} - \sum_{\alpha \in \Delta^+_\theta \setminus \Pi} n_\alpha m_{\alpha,i}$} for $i=1,2,\dots,r$. Since $b_\alpha \in \N_0$ for $\alpha \in \Delta^+$, we get $n_\alpha \in \N_0$ for  $\alpha \in \Delta^+ \setminus \Pi$ and
\begin{align*}
0 \leq a_{\alpha_i}+n_\theta m_{\theta,i} - \sum_{\alpha \in \Delta^+_\theta \setminus \Pi} n_\alpha m_{\alpha,i} \leq a_{\alpha_i}+n_\theta m_{\theta,i}
\end{align*}
for $i=1,2,\dots,r$. Hence, we obtain that $n_\theta \in \N_0$ and \smash{$a+n_\theta t_\theta \in \N_0^{\Delta^+}$}\!. If we denote by $n_a \in \N_0$ the smallest non-negative integer satisfying \smash{$a+n_a t_\theta \in \N_0^{\Delta^+}$}\!, then we may write $n_\theta$ in the form $n_\theta = n_a + n_\theta'$ for $n_\theta' \in \N_0$. This finishes the proof.}

\theorem{\emph{Let $\lambda \in \mfrak{h}^*$. Then we have
\begin{align}
  zv=\chi_{\lambda+\rho}(z)v
\end{align}
for all $z \in \mfrak{Z}(\mfrak{g})$ and $v \in W^\mfrak{g}_\mfrak{b}(\lambda)$.}}

\proof{For $\lambda \in \mfrak{h}^*$ we show that
\begin{align*}
  \pi_{\lambda+\rho}(z)=\chi_{\lambda+\rho}(z)
\end{align*}
for all $z \in \mfrak{Z}(\mfrak{g})$, which implies the required statement. For that reason, let us consider the Verma module $M^\mfrak{g}_\mfrak{b}(\lambda)$ for $\lambda \in \mfrak{h}^*$. Then we have $zv = \chi_{\lambda+\rho}(z)v$ for all $z \in \mfrak{Z}(\mfrak{g})$ and $v \in M^\mfrak{g}_\mfrak{b}(\lambda)$, which follows from \eqref{eq:center action Verma}. As the Verma module $M^\mfrak{g}_\mfrak{b}(\lambda)$ is realized as $\eus{A}^\mfrak{g}_{\widebar{\mfrak{n}}}/\eus{I}_{\rm V} \simeq \C[\partial_{x_\alpha},\, \alpha \in \Delta^+]$, where the $\mfrak{g}$-module structure on $\eus{A}^\mfrak{g}_{\widebar{\mfrak{n}}}/\eus{I}_{\rm V}$ is given through the homomorphism $\pi_{\lambda+\rho} \colon \mfrak{g} \rarr \eus{A}^\mfrak{g}_{\widebar{\mfrak{n}}}$ of Lie algebras, we have
\begin{align*}
  \pi_{\lambda+\rho}(z)v=\chi_{\lambda+\rho}(z)v
\end{align*}
for all $v \in \C[\partial_{x_\alpha},\, \alpha \in \Delta^+]$. Therefore, we have $\pi_{\lambda+\rho}(z)=\chi_{\lambda+\rho}(z)$ for all $z \in \mfrak{Z}(\mfrak{g})$.}

In the following, we give an explicit form of $\pi_\lambda(\Cas_\theta)$ for all $\lambda \in \mfrak{h}^*$. Let us recall that the Lie subalgebra of $\mfrak{g}$ generated by $\{e_\theta,h_\theta,f_\theta\}$ is denoted by $s_\theta$.
\medskip

For $\lambda\in \mfrak{h}^*$ we introduce a linear mapping
\begin{align}
  \sigma_\lambda \colon \mfrak{g} \rarr \eus{A}^\mfrak{g}_{\widebar{\mfrak{n}}}
\end{align}
by
\begin{align}
  \sigma_\lambda(a)= -\sum_{\alpha \in \Delta^+}\bigg[{\ad(v(x))e^{\ad(v(x))} \over e^{\ad(v(x))}-\id}\,(e^{-\ad(v(x))}a)_{\widebar{\mfrak{n}}}\bigg]_\alpha \partial_{x_\alpha} + (\lambda+\rho)((e^{-\ad(v(x))}a)_\mfrak{h}) \label{eq:pi action general res}
\end{align}
for all $a \in \mfrak{g}$, where
\begin{align}
  v(x)= \sum_{\alpha \in \Delta^+_\theta} x_\alpha f_\alpha
\end{align}
and we use the same convention as in \eqref{eq:pi action general}. In particular, we have
\begin{align}
  \sigma_\lambda(a)= - \sum_{\alpha \in\Delta^+} \bigg[{\ad(v(x)) \over e^{\ad(v(x))} - \id}\,a\bigg]_\alpha \partial_{x_\alpha} \label{eq:pi action nilradical res}
\end{align}
for $a \in \widebar{\mfrak{n}}$ and
\begin{align}
  \sigma_\lambda(a)= \sum_{\alpha \in \Delta^+} [\ad(v(x))a]_\alpha \partial_{x_\alpha} + (\lambda+\rho)(a) \label{eq:pi action cartan res}
\end{align}
for $a \in \mfrak{h}$.
\medskip

\proposition{\label{prop:pi action nilradical}
Let $\lambda \in \mfrak{h}^*$. Then we have
\begin{align}
    \pi_\lambda(e_\alpha)&= x_\theta \pi_\lambda([e_\alpha, f_\theta]) + \sigma_\lambda(e_\alpha)
\end{align}
for $\alpha \in \Delta^+_\theta$.}

\proof{Let us assume that $\alpha \in \Delta^+_\theta$. Since $u(x)=v(x)+x_\theta f_\theta$ and $[v(x),x_\theta f_\theta]=0$, we may write
\begin{align*}
e^{-\ad(u(x))}e_\alpha &= e^{-\ad(v(x))-\ad(x_\theta f_\theta)}e_\alpha = e^{-\ad(v(x))} e^{-\ad(x_\theta f_\theta)}e_\alpha \\
&= e^{-\ad(v(x))}(e_\alpha-\ad(x_\theta f_\theta)e_\alpha) =  e^{-\ad(v(x))}(e_\theta + x_\theta [e_\alpha,f_\theta]),
\end{align*}
which implies
\begin{align}
  (e^{-\ad(u(x))}e_\alpha)_\mfrak{h}= (e^{-\ad(v(x))}e_\alpha)_\mfrak{h} \label{eq:projection e_alpha h}
\end{align}
and
\begin{align}
  (e^{-\ad(u(x))}e_\alpha)_{\widebar{\mfrak{n}}} =   (e^{-\ad(v(x))}e_\alpha)_{\widebar{\mfrak{n}}} +x_\theta e^{-\ad(v(x))}[e_\alpha,f_\theta]. \label{eq:projection e_alpha n}
\end{align}
Finally, using \eqref{eq:pi action general} we may write
\begin{align*}
  \pi_\lambda(e_\alpha)&= -\sum_{\alpha \in \Delta^+}\bigg[{\ad(u(x))e^{\ad(u(x))} \over e^{\ad(u(x))}-\id}\,(e^{-\ad(u(x))}e_\alpha)_{\widebar{\mfrak{n}}}\bigg]_\alpha \partial_{x_\alpha} + (\lambda+\rho)((e^{-\ad(u(x))}e_\alpha)_\mfrak{h}) \\
  & = -\sum_{\alpha \in \Delta^+}\bigg[{\ad(v(x))e^{\ad(v(x))} \over e^{\ad(v(x))}-\id}\,(e^{-\ad(u(x))}e_\alpha)_{\widebar{\mfrak{n}}}\bigg]_\alpha \partial_{x_\alpha} + (\lambda+\rho)((e^{-\ad(u(x))}e_\alpha)_\mfrak{h}) \\
  &= -\sum_{\alpha \in \Delta^+}\bigg[{\ad(v(x))e^{\ad(v(x))} \over e^{\ad(v(x))}-\id}\,(e^{-\ad(v(x))}e_\alpha)_{\widebar{\mfrak{n}}}\bigg]_\alpha \partial_{x_\alpha}  + (\lambda+\rho)((e^{-\ad(v(x))}e_\alpha)_\mfrak{h}) \\
  & \quad -\sum_{\alpha \in \Delta^+} x_\theta \bigg[{\ad(v(x))e^{\ad(v(x))} \over e^{\ad(v(x))}-\id}\,e^{-\ad(v(x))}[e_\alpha,f_\theta]\bigg]_\alpha \partial_{x_\alpha},
\end{align*}
were we used $[u(x),a]=[v(x),a]$ for all $a \in \widebar{\mfrak{n}}$ in the second equality and the formulas \eqref{eq:projection e_alpha h} and \eqref{eq:projection e_alpha n} in the last equality. Therefore, we get
\begin{align*}
  \pi_\lambda(e_\alpha)= -\sum_{\alpha \in \Delta^+} x_\theta \bigg[{\ad(v(x)) \over e^{\ad(v(x))}-\id}\,[e_\alpha,f_\theta]\bigg]_\alpha \partial_{x_\alpha} + \sigma_\lambda(e_\alpha) = x_\theta \pi_\lambda([e_\alpha,f_\theta]) + \sigma_\lambda(e_\alpha).
\end{align*}
This finishes the proof.}

\proposition{\label{prop:pi action sl(2,C)-triple}
We have
\begin{align}
  \begin{aligned}
    \pi_\lambda(f_\theta)&=-\partial_{x_\theta}, \\[3mm]
    \pi_\lambda(h_\theta)&= \sum_{\alpha \in \Delta^+} \alpha(h_\theta) x_\alpha \partial_{x_\alpha} + (\lambda+\rho)(h_\theta),\\
    \pi_\lambda(e_\theta)&= x_\theta(\pi_\lambda(h_\theta)-x_\theta \partial_{x_\theta}) + \sigma_\lambda(e_\theta),
  \end{aligned}
\end{align}
for $\lambda \in \mfrak{h}^*$.}

\proof{Since $f_\theta \in \mfrak{z}(\widebar{\mfrak{n}})$, we have $[u(x),f_\theta]=0$. Hence, from \eqref{eq:pi action nilradical} we obtain
\begin{align*}
  \pi_\lambda(f_\theta)= -\sum_{\alpha \in \Delta^+} [f_\theta]_\alpha \partial_{x_\alpha} = -\partial_{x_\theta}.
\end{align*}
Further, we have
\begin{align*}
  [u(x),h_\theta]= \sum_{\alpha \in \Delta^+} x_\alpha [f_\alpha, h_\theta] = \sum_{\alpha \in \Delta^+} \alpha(h_\theta)x_\alpha f_\alpha,
\end{align*}
which together with \eqref{eq:pi action cartan} gives us
\begin{align*}
  \pi_\lambda(h_\theta)= \sum_{\alpha \in \Delta^+} [u(x),h_\theta]_\alpha \partial_{x_\alpha} + (\lambda+\rho)(h_\theta) = \sum_{\alpha \in \Delta^+} \alpha(h_\theta) x_\alpha \partial_{x_\alpha} + (\lambda+\rho)(h_\theta).
\end{align*}
Since $u(x)=v(x)+x_\theta f_\theta$ and $[v(x),x_\theta f_\theta]=0$, we may write
\begin{align*}
e^{-\ad(u(x))}e_\theta &= e^{-\ad(v(x))-\ad(x_\theta f_\theta)}e_\theta = e^{-\ad(v(x))} e^{-\ad(x_\theta f_\theta)}e_\theta \\
&= e^{-\ad(v(x))}(e_\theta-\ad(x_\theta f_\theta)e_\theta + {\textstyle {1 \over 2}}\ad(x_\theta f_\theta)^2 e_\theta) =  e^{-\ad(v(x))}(e_\theta + x_\theta h_\theta -x_\theta^2 f_\theta),
\end{align*}
which implies
\begin{align}
  (e^{-\ad(u(x))}e_\theta)_\mfrak{h}= (e^{-\ad(v(x))}e_\theta)_\mfrak{h} + x_\theta h_\theta \label{eq:projection e_theta h}
\end{align}
and
\begin{align}
  (e^{-\ad(u(x))}e_\theta)_{\widebar{\mfrak{n}}} =   (e^{-\ad(v(x))}e_\theta)_{\widebar{\mfrak{n}}} +x_\theta (e^{-\ad(v(x))}h_\theta)_{\widebar{\mfrak{n}}}-x_\theta^2 f_\theta. \label{eq:projection e_theta n}
\end{align}
Finally, using \eqref{eq:pi action general} we may write
\begin{align*}
  \pi_\lambda(e_\theta)&= -\sum_{\alpha \in \Delta^+}\bigg[{\ad(u(x))e^{\ad(u(x))} \over e^{\ad(u(x))}-\id}\,(e^{-\ad(u(x))}e_\theta)_{\widebar{\mfrak{n}}}\bigg]_\alpha \partial_{x_\alpha} + (\lambda+\rho)((e^{-\ad(u(x))}e_\theta)_\mfrak{h}) \\
  & = -\sum_{\alpha \in \Delta^+}\bigg[{\ad(v(x))e^{\ad(v(x))} \over e^{\ad(v(x))}-\id}\,(e^{-\ad(u(x))}e_\theta)_{\widebar{\mfrak{n}}}\bigg]_\alpha \partial_{x_\alpha} + (\lambda+\rho)((e^{-\ad(u(x))}e_\theta)_\mfrak{h}) \\
  &= -\sum_{\alpha \in \Delta^+}\bigg[{\ad(v(x))e^{\ad(v(x))} \over e^{\ad(v(x))}-\id}\,(e^{-\ad(v(x))}e_\theta)_{\widebar{\mfrak{n}}}\bigg]_\alpha \partial_{x_\alpha}  + (\lambda+\rho)((e^{-\ad(v(x))}e_\theta)_\mfrak{h}) \\
  & \quad -\sum_{\alpha \in \Delta^+} x_\theta \bigg[{\ad(v(x))e^{\ad(v(x))} \over e^{\ad(v(x))}-\id}\,(e^{-\ad(v(x))}h_\theta)_{\widebar{\mfrak{n}}}\bigg]_\alpha \partial_{x_\alpha}  + \sum_{\alpha \in \Delta^+} x_\theta^2 \bigg[{\ad(v(x))e^{\ad(v(x))} \over e^{\ad(v(x))}-\id}\,f_\theta\bigg]_\alpha \partial_{x_\alpha} \\
  & \quad + x_\theta(\lambda+\rho)(h_\theta),
\end{align*}
were we used $[u(x),a]=[v(x),a]$ for all $a \in \widebar{\mfrak{n}}$ in the second equality and the formulas \eqref{eq:projection e_theta h} and \eqref{eq:projection e_theta n} in the last equality. Therefore, we get
\begin{align*}
  \pi_\lambda(e_\theta)&= -\sum_{\alpha \in \Delta^+} x_\theta \bigg[{\ad(v(x))e^{\ad(v(x))} \over e^{\ad(v(x))}-\id}\,(e^{-\ad(v(x))}h_\theta)_{\widebar{\mfrak{n}}}\bigg]_\alpha \partial_{x_\alpha}  + \sum_{\alpha \in \Delta^+} x_\theta^2 [f_\theta]_\alpha \partial_{x_\alpha} \\
  & \quad + x_\theta(\lambda+\rho)(h_\theta) + \sigma_\lambda(e_\theta) \\
  &= -\sum_{\alpha \in \Delta^+} x_\theta \bigg[{\ad(v(x))e^{\ad(v(x))} \over e^{\ad(v(x))}-\id}\,(e^{-\ad(v(x))}h_\theta)_{\widebar{\mfrak{n}}}\bigg]_\alpha \partial_{x_\alpha}  +  x_\theta^2 \partial_{x_\theta} \\
  & \quad + x_\theta(\lambda+\rho)(h_\theta) + \sigma_\lambda(e_\theta).
\end{align*}
Further, since we have $(e^{-\ad(v(x))}h_\theta)_{\widebar{\mfrak{n}}}= (e^{-\ad(v(x))}-\id)h_\theta$, we may write
\begin{align*}
  \pi_\lambda(e_\theta)&= -\sum_{\alpha \in \Delta^+} x_\theta \bigg[{\ad(v(x))e^{\ad(v(x))} \over e^{\ad(v(x))}-\id}\,(e^{-\ad(v(x))}-\id)h_\theta\bigg]_\alpha \partial_{x_\alpha}  +  x_\theta^2 \partial_{x_\theta} \\
  & \quad + x_\theta(\lambda+\rho)(h_\theta) + \sigma_\lambda(e_\theta) \\
  & = \sum_{\alpha \in \Delta^+} x_\theta [\ad(v(x)) h_\theta]_\alpha \partial_{x_\alpha}  +  x_\theta^2 \partial_{x_\theta} + x_\theta(\lambda+\rho)(h_\theta) + \sigma_\lambda(e_\theta) \\
  &=\sum_{\alpha \in \Delta^+} x_\theta [u(x),h_\theta]_\alpha \partial_{x_\alpha}  -  x_\theta^2 \partial_{x_\theta} + x_\theta(\lambda+\rho)(h_\theta) + \sigma_\lambda(e_\theta) \\
  &= x_\theta(\pi_\lambda(h_\theta)-x_\theta \partial_{x_\theta}) + \sigma_\lambda(e_\theta).
\end{align*}
This finishes the proof.}

\theorem{\label{thm:cyclic}
\emph{Let $\lambda \in \mfrak{h}^*$. Then $W^\mfrak{g}_\mfrak{b}(\lambda)$ is a cyclic $\mfrak{g}$-module.}}

\proof{Since the polynomials $[(e^{-\ad(v(x))}e_\theta)_\mfrak{h}]_\alpha$ and $[(e^{-\ad(v(x))}e_\theta)_{\widebar{\mfrak{n}}}]_\alpha$ do not contain the constant and linear term for all $\alpha \in \Delta^+$, we obtain, using \eqref{eq:pi action general res}, that $\sigma_{\lambda+\rho}(e_\theta)x_\theta^n=0$ for $n \in \N_0$. Hence, we have
\begin{align*}
  \pi_{\lambda+\rho}(e_\theta)x_\theta^n &= x_\theta(\pi_{\lambda+\rho}(h_\theta) - x_\theta \partial_{x_\theta})x_\theta^n = x_\theta(x_\theta \partial_{x_\theta} +(\lambda + \theta)(h_\theta)) x_\theta^n = (\lambda(h_\theta)+ n +2) x_\theta^{n+1}, \\
  \pi_{\lambda+\rho}(h_\theta)x_\theta^n &= (2x_\theta \partial_{x_\theta}+(\lambda+\theta)(h_\theta))x_\theta^n = (\lambda(h_\theta)+2+2n)x_\theta^n, \\
  \pi_{\lambda+\rho}(f_\theta)x_\theta^n &= -n x_\theta^{n-1}
\end{align*}
for $n \in \N_0$, where we used \eqref{eq:action weight}, which gives us $\C[x_\theta]=U(\mfrak{s}_\theta)1$ provided $\lambda(h_\theta)+2 \notin -\N_0$ and
\smash{$\C[x_\theta]=U(\mfrak{s}_\theta)x_\theta^{-\lambda(h_\theta)-1}$} if
$\lambda(h_\theta)+2 \in -\N_0$.

From \eqref{eq:pi action nilradical} and the identity
\begin{align*}
  {\ad(u(x)) \over e^{\ad(u(x))}-\id}\,a = a + {\ad(u(x))-e^{\ad(u(x))} + \id \over e^{\ad(u(x))}-\id}\,a
\end{align*}
for all $a \in \widebar{\mfrak{n}}$, we obtain
\begin{align*}
  \pi_{\lambda+\rho}(f_\alpha) = -\partial_{x_\alpha} + \pi'_{\lambda+\rho}(f_\alpha)
\end{align*}
for $\alpha \in \Delta^+$, where
\begin{align*}
  \pi'_{\lambda+\rho}(a)= -\sum_{\alpha\in\Delta^+} \bigg[{\ad(u(x))-e^{\ad(u(x))} + \id \over e^{\ad(u(x))}-\id}\,a \bigg]_\alpha \partial_{x_\alpha}.
\end{align*}
Further, for $k \in \N_0$ we denote by $F_k W^\mfrak{g}_\mfrak{b}(\lambda)$ the vector subspace of $W^\mfrak{g}_\mfrak{b}(\lambda)$ consisting of all polynomials of degree at most $k$ in the variables $\partial_{x_\alpha}$, $\alpha \in \Delta^+_\theta$. Hence, we have an increasing filtration $\{F_k W^\mfrak{g}_\mfrak{b}(\lambda)\}_{k \in \N_0}$ on $W^\mfrak{g}_\mfrak{b}(\lambda)$. As the polynomial
\begin{align*}
  \bigg[{\ad(u(x)) - e^{\ad(u(x))} + \id \over e^{\ad(u(x))} - \id}\,a\bigg]_\alpha
\end{align*}
does not contain the constant term and does not depend on $x_\theta$ for all $\alpha \in \Delta^+$ and $a \in \widebar{\mfrak{n}}$, we get
\begin{align*}
  \pi'_{\lambda+\rho}(f_\alpha)(F_k W^\mfrak{g}_\mfrak{b}(\lambda)) \subset F_k W^\mfrak{g}_\mfrak{b}(\lambda)
\end{align*}
for $k \in \N_0$ and $\alpha \in \Delta^+$. Therefore, we have
\begin{align*}
  \pi_{\lambda+\rho}(f_\alpha)(F_k W^\mfrak{g}_\mfrak{b}(\lambda))/F_k W^\mfrak{g}_\mfrak{b}(\lambda)= \partial_{x_\alpha} F_k W^\mfrak{g}_\mfrak{b}(\lambda)/ F_k W^\mfrak{g}_\mfrak{b}(\lambda)
\end{align*}
for $k \in \N_0$ and $\alpha \in \Delta^+$, which together with $F_0 W^\mfrak{g}_\mfrak{b}(\lambda)= \C[x_\theta]$ gives us that $W^\mfrak{g}_\mfrak{b}(\lambda)$ is generated by $\C[x_\theta]$. As $\C[x_\theta]$ is generated by one vector, we get that $W^\mfrak{g}_\mfrak{b}(\lambda)$ is a cyclic $\mfrak{g}$-module.}

\proposition{\label{prop:Casimir description}
Let $\lambda \in \mfrak{h}^*$. Then the Casimir operator $\pi_\lambda(\Cas_\theta)$ has the form
\begin{align}
  \pi_\lambda(\Cas_\theta)= \pi_\lambda(\Cas_\theta)_s + \pi_\lambda(\Cas_\theta)_n,
\end{align}
where
\begin{align}
  \pi_\lambda(\Cas_\theta)_s= {\textstyle {1 \over 2}} \sigma_\lambda(h_\theta)(\sigma_\lambda(h_\theta)-2)
\end{align}
is a semisimple operator and
\begin{align}
  \pi_\lambda(\Cas_\theta)_n= 2\sigma_\lambda(e_\theta)\sigma_\lambda(f_\theta)
\end{align}
is a locally nilpotent operator. Moreover, the weight vector $w_{\lambda,a+ (n_a+n_\theta)t_\theta+t}$ with weight $\mu_{\lambda,a}$ for $a \in \Z^\Pi$, $n_\theta \in \N_0$ and $t \in \Lambda_+^\theta$ is an eigenvector of \smash{$\pi_{\lambda+\rho}(\Cas_\theta)_s$} with eigenvalue
\begin{align}
  c_{\lambda,a,n_\theta} = {\textstyle {1 \over 2}} (\mu_{\lambda,a}(h_\theta)- 2(n_a+n_\theta))(\mu_{\lambda,a}(h_\theta)- 2(n_a+n_\theta+1))
\end{align}
provided \smash{$a+(n_a+n_\theta)t_\theta+t \in \N_0^{\Delta^+}$}\!, where $n_a \in \N_0$ is the smallest non-negative integer satisfying \smash{$a+n_at_\theta \in \N_0^{\Delta^+}$}\!.}

\proof{As $\Cas_\theta=e_\theta f_\theta + f_\theta e_\theta + {1 \over 2} h_\theta^2 = 2f_\theta e_\theta + {1 \over 2} h_\theta(h_\theta +2)$, we may write
\begin{align*}
  \pi_\lambda(\Cas_\theta) &= 2 \pi_\lambda(f_\theta) \pi_\lambda(e_\theta) + {\textstyle {1 \over 2}} \pi_\lambda(h_\theta)(\pi_\lambda(h_\theta)+2) \\
  &= -2\partial_{x_\theta}(x_\theta(\pi_\lambda(h_\theta)-x_\theta \partial_{x_\theta})+ \sigma_\lambda(e_\theta)) + {\textstyle {1 \over 2}} \pi_\lambda(h_\theta)(\pi_\lambda(h_\theta)+2) \\
  &= -2(x_\theta \partial_{x_\theta}+1)(\pi_\lambda(h_\theta)-x_\theta \partial_{x_\theta}) - 2\partial_{x_\theta} \sigma_\lambda(e_\theta) + {\textstyle {1 \over 2}} \pi_\lambda(h_\theta)(\pi_\lambda(h_\theta)+2),
\end{align*}
where we used Proposition \ref{prop:pi action sl(2,C)-triple}. If we denote
\begin{align*}
  \pi_\lambda(\Cas_\theta)_s = -2(x_\theta \partial_{x_\theta}+1)(\pi_\lambda(h_\theta)-x_\theta \partial_{x_\theta}) + {\textstyle {1 \over 2}} \pi_\lambda(h_\theta)(\pi_\lambda(h_\theta)+2)
\end{align*}
and
\begin{align*}
  \pi_\lambda(\Cas_\theta)_n = -2 \partial_{x_\theta} \sigma_\lambda(e_\theta),
\end{align*}
we obtain $\pi_\lambda(\Cas_\theta) = \pi_\lambda(\Cas_\theta)_s + \pi_\lambda(\Cas_\theta)_n$. Further, we have
\begin{align*}
  \pi_\lambda(\Cas_\theta)_n = -2 \partial_{x_\theta} \sigma_\lambda(e_\theta) = -2\sigma_\lambda(e_\theta) \partial_{x_\theta} = 2\sigma_\lambda(e_\theta)\sigma_\lambda(f_\theta)
\end{align*}
and
\begin{align*}
\pi_\lambda(\Cas_\theta)_s &= -2(x_\theta \partial_{x_\theta}+1)(\pi_\lambda(h_\theta)-x_\theta \partial_{x_\theta}) + {\textstyle {1 \over 2}} \pi_\lambda(h_\theta)(\pi_\lambda(h_\theta)+2) \\
&= -2(x_\theta \partial_{x_\theta}+1)(\sigma_\lambda(h_\theta)+x_\theta \partial_{x_\theta}) + {\textstyle {1 \over 2}} (\sigma_\lambda(h_\theta) + 2 x_\theta \partial_{x_\theta})(\sigma_\lambda(h_\theta) + 2 x_\theta \partial_{x_\theta}+2) \\
&= -2\sigma_\lambda(h_\theta) + {\textstyle {1 \over 2}} \sigma_\lambda(h_\theta)(\sigma_\lambda(h_\theta+2) = {\textstyle {1 \over 2}} \sigma_\lambda(h_\theta)(\sigma_\lambda(h_\theta)-2),
\end{align*}
where we used $\pi_\lambda(h_\theta) = \sigma_\lambda(h_\theta) + 2 x_\theta \partial_{x_\theta}$, which easily follows from \eqref{eq:pi action cartan res} and Proposition \ref{prop:pi action sl(2,C)-triple}.

Now, let $a \in \Z^\Pi$. Then the weight vector $w_{\lambda-\rho,a+(n_a+n_\theta)t_\theta+ t}$ with weight $\mu_{\lambda-\rho,a}$ for $n_\theta \in \N_0$ and $t \in \Lambda_+^\theta$ is an eigenvector of $\sigma_\lambda(h_\theta)$ and therefore also of $\pi_\lambda(\Cas_\theta)_s$ if \smash{$a+(n_a+n_\theta)t_\theta+t \in \N_0^{\Delta^+}$}\!, since we have
\begin{align*}
  \sigma_\lambda(h_\theta)w_{\lambda-\rho,a+(n_a+n_\theta)t_\theta + t} &= (\pi_\lambda(h_\theta)-2x_\theta \partial_{x_\theta})w_{\lambda-\rho,a+(n_a+n_\theta)t_\theta + t} \\
  &=(\mu_{\lambda-\rho,a}(h_\theta)- 2(n_a+n_\theta))w_{\lambda-\rho,a+(n_a+n_\theta)t_\theta + t}.
\end{align*}
As $\{w_{\lambda-\rho,a+(n_a+n_\theta)t_\theta+t};\, a \in \Z^\Pi,\, n_\theta \in \N_0,\, t \in \Lambda_+^\theta,\, a+(n_a+n_\theta)t_\theta + t \in \smash{\N_0^{\Delta^+}}\!\}$ forms a basis of $W^\mfrak{g}_\mfrak{b}(\lambda-\rho)$, we get that $\pi_\lambda(\Cas_\theta)_s$ is a semisimple operator.

Since $\sigma_\lambda(f_\theta)$ is obviously a locally nilpotent operator and $\sigma_\lambda(f_\theta)\sigma_\lambda(e_\theta)= \sigma_\lambda(e_\theta) \sigma_\lambda(f_\theta)$, we obtain that $\pi_\lambda(\Cas_\theta)_n$ is also a locally nilpotent operator.}

The main disadvantage of $\mfrak{g}$-modules $W^\mfrak{g}_\mfrak{b}(\lambda)$ for $\lambda \in \mfrak{h}^*$, in comparison with objects of the category $\mcal{O}$, is infinite dimensional weight spaces. However, we can enlarge the Cartan subalgebra $\mfrak{h}$ of $\mfrak{g}$ to a commutative subalgebra $\Gamma_\theta$ of $U(\mfrak{g})$  in such a way that $W^\mfrak{g}_\mfrak{b}(\lambda)$ for $\lambda \in \mfrak{h}^*$ is a  $\theta$-Gelfand--Tsetlin module, i.e.\ $W^\mfrak{g}_\mfrak{b}(\lambda) \in \mcal{H}(\mathfrak{g}, \Gamma_\theta)$.
\medskip

\theorem{\emph{Let $\lambda \in \mfrak{h}^*$. Then the $\mfrak{g}$-module $W^\mfrak{g}_\mfrak{b}(\lambda)$ is a  $\theta$-Gelfand--Tsetlin module with finite $\Gamma_\theta$-multiplicities, i.e.\ $W^\mfrak{g}_\mfrak{b}(\lambda) \in \mcal{H}_{\rm fin}(\mfrak{g}, \Gamma_\theta)$.}}

\proof{To prove that $\Gamma_\theta$-weight spaces are finite-dimensional we apply Theorem \ref{thm:generalized eigenspace isomorphism}. For that reason, let us introduce an $\N_0$-grading on $W^\mfrak{g}_\mfrak{b}(\lambda)$. For $r \in \N_0$ we denote by $W^\mfrak{g}_\mfrak{b}(\lambda)_r$ the vector subspace of $W^\mfrak{g}_\mfrak{b}(\lambda)$ consisting of all homogenous polynomials of degree $r$ in the variables $\partial_{x_\alpha}$, $\alpha \in \Delta^+_\theta$, and $x_\theta$. Moreover, we have $\dim W^\mfrak{g}_\mfrak{b}(\lambda)_r < \infty$ for all $r \in \N_0$. Since the polynomials $[(e^{-\ad(v(x))}e_\theta)_\mfrak{h}]_\alpha$ and $[(e^{-\ad(v(x))}e_\theta)_{\widebar{\mfrak{n}}}]_\alpha$ do not contain the constant and linear term for all $\alpha \in \Delta^+$, we obtain, using \eqref{eq:pi action general res}, that
\begin{align*}
   \sigma_{\lambda+\rho}(e_\theta)(F_r W^\mfrak{g}_\mfrak{b}(\lambda)) \subset F_{r-1}W^\mfrak{g}_\mfrak{b}(\lambda)
\end{align*}
for $r \in \N_0$, which together with
\begin{align*}
  \sigma_{\lambda+\rho}(h_\theta)(W^\mfrak{g}_\mfrak{b}(\lambda)_r) \subset W^\mfrak{g}_\mfrak{b}(\lambda)_r \qquad \text{and} \qquad \sigma_{\lambda+\rho}(f_\theta)(W^\mfrak{g}_\mfrak{b}(\lambda)_r) \subset W^\mfrak{g}_\mfrak{b}(\lambda)_{r-1}
\end{align*}
for $r \in \N_0$, according to Proposition \ref{prop:Casimir description}, gives us
\begin{align*}
  \pi_{\lambda+\rho}(\Cas_\theta)_s(W^\mfrak{g}_\mfrak{b}(\lambda)_r) \subset W^\mfrak{g}_\mfrak{b}(\lambda)_r, \qquad
  \pi_{\lambda+\rho}(\Cas_\theta)_s(F_rW^\mfrak{g}_\mfrak{b}(\lambda)) \subset F_{r-2}W^\mfrak{g}_\mfrak{b}(\lambda)
\end{align*}
for all $r \in \N_0$. Therefore, the assumptions of Theorem \ref{thm:generalized eigenspace isomorphism} are satisfied.

Let $\chi \in \Hom(\Gamma_\theta,\C)$ be a $\Gamma_\theta$-weight of $W^\mfrak{g}_\mfrak{b}(\lambda)$. Let us denote $\mu=\chi_{|\mfrak{h}}$ and $c= \chi(\Cas_\theta)$. Then $\mu \in \mfrak{h}^*$ is an $\mfrak{h}$-weight of $W^\mfrak{g}_\mfrak{b}(\lambda)$ and we have $W^\mfrak{g}_\mfrak{b}(\lambda)_\chi \subset W^\mfrak{g}_\mfrak{b}(\lambda)_\mu$. By Proposition \ref{prop:weights description} we obtain that there exists $a \in \Z^\Pi$ such that $\mu=\mu_{\lambda,a}$. The $\N_0$-grading on $W^\mfrak{g}_\mfrak{b}(\lambda)$ induces the $\N_0$-grading on $W^\mfrak{g}_\mfrak{b}(\lambda)_{\mu_a}$ defined by $W^\mfrak{g}_\mfrak{b}(\lambda)_{\mu_a,r} = W^\mfrak{g}_\mfrak{b}(\lambda)_{\mu_a} \cap W^\mfrak{g}_\mfrak{b}(\lambda)_r$ for all $r \in \N_0$. Hence, by using Theorem \ref{thm:generalized eigenspace isomorphism} for the linear mapping $\pi_{\lambda+\rho}(\Cas_\theta)$ on the vector subspace $W^\mfrak{g}_\mfrak{b}(\lambda)_{\mu_a}$, we get that eigenvalues of $\pi_{\lambda+\rho}(\Cas_\theta)$ on $W^\mfrak{g}_\mfrak{b}(\lambda)_{\mu_a}$ are the same as eigenvalues of $\pi_{\lambda+\rho}(\Cas_\theta)_s$ on $W^\mfrak{g}_\mfrak{b}(\lambda)_{\mu_a}$, and moreover the dimensions of generalized eigenspaces corresponding to the same eigenvalue are equal.
Hence, by Proposition \ref{prop:Casimir description} there exists $n_\theta \in \N_0$ such that $c=c_{\lambda,a,n_\theta}$. Further, from Proposition \ref{prop:weights description} and Proposition \ref{prop:Casimir description} follows that the set
\begin{align*}
  \{w_{\lambda,a+(n_a+n_\theta)t_\theta+t};\, n_\theta \in \N_0,\, t \in \Lambda_+^\theta,\, a+(n_a+n_\theta)t_\theta+t \in \N_0^{\Delta^+}\!\}
\end{align*}
forms a basis of $W^\mfrak{g}_\mfrak{b}(\lambda)_{\mu_a}$ and $w_{\lambda,a+(n_a+n_\theta)t_\theta+t}$ is an eigenvectors of $\pi_{\lambda+\rho}(\Cas_\theta)_s$ with eigenvalue $c_{\lambda,a,n_\theta}$. As we have $c_{\lambda,a,n_\theta} = c_{\lambda,a,n'_\theta}$ for $n_\theta, n'_\theta \in \N_0$ if and only if $n'_\theta =n_\theta$ or $n'_\theta=\mu_{\lambda,a}(h_\theta)-2n_a-n_\theta-1$ provided $\mu_{\lambda,a}(h_\theta)-2n_a-n_\theta-1 \in \N_0$, we obtain
\begin{multline*}
  \dim W^\mfrak{g}_\mfrak{b}(\lambda)_\chi = \natural \{a+(n_a+n_\theta)t_\theta + t;\, t \in \Lambda_+^\theta,\ a+(n_a+n_\theta)t_\theta +t \in \N_0^{\Delta^+}\!\} \\
   + \natural \{a+(\mu_{\lambda,a}(h_\theta)-n_a-n_\theta-1)t_\theta+t;\, t \in \Lambda_+^\theta,\, a+(\mu_{\lambda,a}(h_\theta)-n_a-n_\theta-1)t_\theta+t \in \N_0^{\Delta^+}\!\},
\end{multline*}
which is obviously a finite number. This completes the proof.}

\corollary{\label{cor-clas}
The $\mfrak{g}$-module $W^\mfrak{g}_\mfrak{b}(\lambda)$ belongs to the category $\mcal{I}_{\theta, f, {\rm fin}}$ for all $\lambda \in \mfrak{h}^*$.}

\corollary{\label{cor-GT}
Let $\Gamma$ be any commutative subalgebra of $U(\mfrak{g})$ containing $\Gamma_{\theta}$. Then
the $\mfrak{g}$-module $W^\mfrak{g}_\mfrak{b}(\lambda)$ belongs to the category $\mcal{H}_{\rm fin}(U(\mfrak{g}), \Gamma) $ for all $\lambda \in \mfrak{h}^*$.}

Hence, the $\theta$-Gelfand--Tsetlin $\mfrak{g}$-modules $W^\mfrak{g}_\mfrak{b}(\lambda)$ is  a $\Gamma$-weight module for any commutative $\Gamma$ containing $\Gamma_{\theta}$. In particular, if $\mfrak{g}=\mfrak{sl}(n)$ (or $\mfrak{g}=\mfrak{gl}(n)$) and $\Gamma$ is a Gelfand-Tsetlin subalgebra of $U(\mfrak{g})$ containing $\Cas_\theta$
 then $W^\mfrak{g}_\mfrak{b}(\lambda)$ is a Gelfand-Tsetlin module (with respect to a particular choice of a Gelfand-Tsetlin subalgebra) \cite{Drozd-Futorny-Ovsienko1994} (in this case the chain of $\mfrak{sl}(m)$-subalgebras contains $\mfrak{sl}(2)$ associated to the root $\theta$). Therefore, the $\mfrak{g}$-modules
 $W^\mfrak{g}_\mfrak{b}(\lambda)$ can be viewed as examples of Gelfand-Tsetlin modules for an arbitrary simple finite-dimensional $\mfrak{g}$ with respect to any $\Gamma$ containing $\Gamma_{\theta}$.

 By duality we can similarly construct certain $\theta$-Gelfand--Tsetlin  modules in the category $\mcal{I}_{\theta, e, {\rm fin}}$. A comparison of the basic characteristics of Verma modules and Gelfand--Tsetlin modules is given in Table \ref{tab:comparision}. Their realization is given by the following theorem.
\medskip

\theorem{\label{thm-main}
\emph{Let $\lambda \in \mfrak{h}^*$. Then the $\mfrak{g}$-module
\begin{align}
 \eus{A}^\mfrak{g}_{\widebar{\mfrak{n}}}/ \eus{I}_{\rm GT} \simeq
\C[\partial_{x_\alpha},\alpha \in \Delta^+_\theta,\, x_\theta], \quad \eus{I}_{\rm GT}=(x_\alpha, \alpha \in \Delta^+_\theta,\partial_{x_\theta}),
\end{align}
where the $\mfrak{g}$-module structure on $\eus{A}^\mfrak{g}_{\widebar{\mfrak{n}}}/ \eus{I}_{\rm GT}$ is given through the homomorphism
$\pi_{\lambda+\rho} \colon  \mfrak{g}\rarr \eus{A}^\mfrak{g}_{\widebar{\mfrak{n}}}$ of Lie algebras with
\begin{align}
\pi_\lambda(a)= -\sum_{\alpha \in \Delta^+}\bigg[{\ad(u(x))e^{\ad(u(x))} \over e^{\ad(u(x))}-\id}\,(e^{-\ad(u(x))}a)_{\widebar{\mfrak{n}}}\bigg]_\alpha \partial_{x_\alpha} + (\lambda+\rho)((e^{-\ad(u(x))}a)_\mfrak{h})
\end{align}
for all $a \in \mfrak{g}$, is a weight, cyclic $\theta$-Gelfand--Tsetlin module with finite $\Gamma_\theta$-multiplicities.}}
\vspace{-2mm}


\subsection{Gelfand--Tsetlin modules for $\mfrak{sl}(3,\C)$}

In this section, we will focus more closely on the Gelfand--Tsetlin modules $W^\mfrak{g}_\mfrak{b}(\lambda)$ for $\lambda \in \mfrak{h}^*$ over the Lie algebra $\mfrak{sl}(3,\C)$. In particular, we describe the irreducibility condition for these modules.
We use the notation introduced in the previous section.
\medskip

Let us consider the complex simple Lie algebra $\mfrak{g}=\mfrak{sl}(3,\C)$. A Cartan subalgebra $\mfrak{h}$ of $\mfrak{g}$ is given by diagonal matrices
\begin{align}
  \mfrak{h}=\{\diag(a_1,a_2,a_3);\, a_1,a_2,a_3 \in \C,\ a_1+a_2+a_3=0\}.
\end{align}
For $i=1,2,3$ we define $\veps_i \in \mfrak{h}^*$ by $\veps_i(\diag(a_1,a_2,a_3))=a_i$. The root system of $\mfrak{g}$ with respect to $\mfrak{h}$ is $\Delta=\{\veps_i-\veps_j;\, 1\leq i \neq j \leq 3\}$. A positive root system in $\Delta$ is $\Delta^+=\{\veps_1 - \veps_2, \veps_2-\veps_3, \veps_1-\veps_3\}$ and the set of simple roots is $\Pi=\{\alpha_1,\alpha_2\}$ with $\alpha_1= \veps_1 - \veps_2$ and $\alpha_2=\veps_2-\veps_3$. The remaining positive root $\alpha_1+\alpha_2= \veps_1-\veps_3$ is the maximal root $\theta$ of $\mfrak{g}$.

Let us introduce a root basis of the niradical $\mfrak{n}$ by
\begin{align*}
  e_{\alpha_1}=
  \begin{pmatrix}
    0 & 1 & 0 \\
    0 & 0 & 0 \\
    0 & 0 & 0
  \end{pmatrix}\!, \qquad \qquad
  e_\theta&=
  \begin{pmatrix}
    0 & 0 & 1  \\
    0 & 0 & 0 \\
    0 & 0 & 0
  \end{pmatrix}\!, \qquad \qquad
  e_{\alpha_2}=
  \begin{pmatrix}
    0 & 0 & 0  \\
    0 & 0 & 1 \\
    0 & 0 & 0
  \end{pmatrix}
\end{align*}
and a root basis of the opposite nilradical $\widebar{\mfrak{n}}$ through
\begin{align*}
  f_{\alpha_1}=
  \begin{pmatrix}
    0 & 0 & 0 \\
    1 & 0 & 0 \\
    0 & 0 & 0
  \end{pmatrix}\!, \qquad \qquad
  f_\theta=
  \begin{pmatrix}
    0 & 0 & 0  \\
    0 & 0 & 0 \\
    1 & 0 & 0
  \end{pmatrix}\!, \qquad \qquad
  f_{\alpha_2}=
  \begin{pmatrix}
    0 & 0 & 0  \\
    0 & 0 & 0 \\
    0 & 1 & 0
  \end{pmatrix}\!.
\end{align*}
The coroots corresponding to the positive roots are given by
\begin{align*}
  h_{\alpha_1}=
  \begin{pmatrix}
    1 & 0 & 0 \\
    0 & -1 & 0 \\
    0 & 0 & 0
  \end{pmatrix}\!, \qquad \qquad
  h_\theta=
  \begin{pmatrix}
    1 & 0 & 0 \\
    0 & 0 & 0 \\
    0 & 0 & -1
  \end{pmatrix}\!, \qquad \qquad
  h_{\alpha_2}=
  \begin{pmatrix}
    0 & 0 & 0  \\
    0 & 1 & 0 \\
    0 & 0 & -1
  \end{pmatrix}\!.
\end{align*}
These generators fulfill, among others, the commutation relations $[e_{\alpha_1}, e_{\alpha_2}]=e_\theta$ and $[f_{\alpha_2}, f_{\alpha_1}]=f_\theta$. Moreover, we have that $(e_\theta,h_\theta,f_\theta)$ is an $\mfrak{sl}(2,\C)$-triple.

Let us recall that $\{x_{\alpha_1},x_{\alpha_2},x_\theta\}$ are the linear coordinate functions on $\widebar{\mfrak{n}}$ with respect to the basis $\{f_{\alpha_1}, f_{\alpha_2}, f_\theta\}$ of $\widebar{\mfrak{n}}$ and that the Weyl algebra $\eus{A}^\mfrak{g}_{\widebar{\mfrak{n}}}$ of $\widebar{\mfrak{n}}$ is generated by  $\{x_{\alpha_1},\partial_{x_{\alpha_1}},x_{\alpha_2},\partial_{x_{\alpha_2}},x_\theta,\partial_{x_\theta}\}$ together with the canonical commutation relations. For clarity, we denote $x=x_{\alpha_1}$, $y=x_{\alpha_2}$ and $z=x_\theta$.
\medskip

\theorem{\label{thm:embedding sl(3,C)}
\emph{Let $\lambda \in \mfrak{h}^*$. Then the homomorphism $\pi_\lambda \colon \mfrak{g} \rarr \eus{A}^\mfrak{g}_{\widebar{\mfrak{n}}}$ of Lie algebras is given by
\begin{align*}
\pi_\lambda(f_{\alpha_1})&=-\partial_x+{\textstyle {1 \over 2}}y\partial_z, \\
\pi_\lambda(f_{\alpha_2})&=-\partial_y-{\textstyle {1 \over 2}}x\partial_z, \\
\pi_\lambda(f_\theta)&=-\partial_z, \\
\pi_\lambda(e_{\alpha_1})&= z(\partial_y  + {\textstyle {1 \over 2}}x\partial_z) +  x(x\partial_x - {\textstyle {1\over 2}} y\partial_y + \lambda_1+1)+ {\textstyle {1\over 4}}x^2y\partial_z, \\
\pi_\lambda(e_{\alpha_2})&=z(-\partial_x + {\textstyle {1 \over 2}} y\partial_z) + y(y\partial_y - {\textstyle {1 \over 2}} x \partial_x +  \lambda_2+1) - {\textstyle {1\over 4}}xy^2\partial_z, \\
\pi_\lambda(e_\theta)&=z(x\partial_x + y\partial_y +z\partial_z + \lambda_1+\lambda_2+2) + {\textstyle {1 \over 2}}xy(x\partial_x-y\partial_y+\lambda_1-\lambda_2) +
  {\textstyle {1 \over 4}}x^2y^2\partial_z, \\
\pi_\lambda(h_{\alpha_1})&=2x\partial_x-y\partial_y+z\partial_z+\lambda_1+1, \\
\pi_\lambda(h_{\alpha_2})&=-x\partial_x+2y\partial_y+z\partial_z+\lambda_2+1, \\
\pi_\lambda(h_\theta)&= x\partial_x+y\partial_y+2z\partial_z+\lambda_1+ \lambda_2+2,
\end{align*}
where $\lambda_1=\lambda(h_{\alpha_1})$ and $\lambda_2=\lambda(h_{\alpha_2})$.}}

\proof{It follows immediately by a straightforward computation from \eqref{eq:pi action general}.}

The homomorphism $\pi_\lambda \colon \mfrak{g} \rarr \eus{A}^\mfrak{g}_{\widebar{\mfrak{n}}}$ of Lie algebras enables us to construct a wide class of interesting $\mfrak{g}$-modules. The easiest way is to take a left ideal $\eus{I}$ of $\eus{A}^\mfrak{g}_{\widebar{\mfrak{n}}}$. Then $\eus{A}^\mfrak{g}_{\widebar{\mfrak{n}}}/\eus{I}$ has the $\mfrak{g}$-module structure through the homomorphism $\pi_\lambda \colon \mfrak{g} \rarr \eus{A}^\mfrak{g}_{\widebar{\mfrak{n}}}$.

In the previous section, we considered only the two simplest possible left ideals of the Weyl algebra $\eus{A}^\mfrak{g}_{\widebar{\mfrak{n}}}$. The first one $\eus{I}_{\rm V}=(x,y,z)$ leads to the Verma modules $M^\mfrak{g}_\mfrak{b}(\lambda) \simeq \eus{A}^\mfrak{g}_{\widebar{\mfrak{n}}}/\eus{I}_{\rm V}$, while the second one $\eus{I}_{\rm GT}=(x,y,\partial_z)$ leads to the Gelfand--Tsetlin modules $W^\mfrak{g}_\mfrak{b}(\lambda) \simeq \eus{A}^\mfrak{g}_{\widebar{\mfrak{n}}}/\eus{I}_{\rm GT}$. However, we have much more options, where some of them lead to the so called twisted Verma modules $M^\mfrak{g}_\mfrak{b}(\lambda)^w$ labeled by the elements $w \in W$ of the Weyl group $W$ of $\mfrak{g}$, see \cite{Krizka-Somberg2015b}.

The twisted Verma modules $M^\mfrak{g}_\mfrak{b}(\lambda)^w$ for $\lambda \in \mfrak{h}^*$ and $w \in W$ are realized as $\eus{A}^\mfrak{g}_{\widebar{\mfrak{n}}}/\eus{I}_w$, where
$\eus{I}_w$ is the left ideal of $\eus{A}^\mfrak{g}_{\widebar{\mfrak{n}}}$ defined by
\begin{align}
\eus{I}_w = (x_\alpha,\,\alpha\in w^{-1}(\Delta^+)\cap\Delta^+,\ \partial_{x_\alpha},\,\alpha\in w^{-1}(-\Delta^+)\cap\Delta^+)
\end{align}
and the $\mfrak{g}$-module structure on $\eus{A}^\mfrak{g}_{\widebar{\mfrak{n}}}/\eus{I}_w$ is given through the homomorphism $\pi^w_{\lambda+\rho} \colon \mfrak{g} \rarr \eus{A}^\mfrak{g}_{\widebar{\mfrak{n}}}$ of Lie algebras defined by
\begin{align}
  \pi^w_{\lambda+\rho}=\pi_{w^{-1}(\lambda+\rho)}\circ \Ad(\dot{w}^{-1}),
\end{align}
where $w^{-1}$ acts in the standard way and $\dot{w}$ is a representative of $w$ in the Lie group $G=\SL(3,\C)$. The list of all possibilities looks as follows:
\begin{enumerate}
  \item[1)] $w=e$, $\eus{I}_e=(x,y,z)$, $\eus{A}^\mfrak{g}_{\widebar{\mfrak{n}}}/I_e \simeq
\C[\partial_x, \partial_y, \partial_z]$;
  \item[2)] $w=s_1$,  $\eus{I}_{s_1}=(\partial_x,y,z)$, $\eus{A}^\mfrak{g}_{\widebar{\mfrak{n}}}/\eus{I}_{s_1}
\simeq \C[x, \partial_y, \partial_z]$;
  \item[3)] $w=s_2$, $\eus{I}_{s_2}=(x,\partial_y,z)$, $\eus{A}^\mfrak{g}_{\widebar{\mfrak{n}}}/\eus{I}_{s_2}
\simeq \C[\partial_x, y, \partial_z]$;
  \item[4)] $w=s_1s_2$, $\eus{I}_{s_1s_2}=(x,\partial_y,\partial_z)$, $\eus{A}^\mfrak{g}_{\widebar{\mfrak{n}}}/\eus{I}_{s_1s_2}
\simeq \C[\partial_x, y, z]$;
  \item[5)] $w=s_2s_1$, $\eus{I}_{s_2s_1}=(\partial_x,y,\partial_z)$, $\eus{A}^\mfrak{g}_{\widebar{\mfrak{n}}}/\eus{I}_{s_2s_1}
\simeq \C[x,\partial_y, z]$;
  \item[6)] $w=s_1s_2s_1$, $\eus{I}_{s_1s_2s_1}=(\partial_x,\partial_y,\partial_z)$, $\eus{A}^\mfrak{g}_{\widebar{\mfrak{n}}}/\eus{I}_{s_1s_2s_1}
\simeq \C[x, y, z]$.
\end{enumerate}
The elements $s_1$ and $s_2$ of the Weyl group $W$ denote the reflections about the hyperplanes perpendicular to $\alpha_1$ and $\alpha_2$, respectively.

In fact, we have another two possible left ideals of $\eus{A}^\mfrak{g}_{\widebar{\mfrak{n}}}$, which look as follows:
\begin{enumerate}
  \item[7)] $\eus{I}_{\rm GT}=(x,y,\partial_z)$, $\eus{A}^\mfrak{g}_{\widebar{\mfrak{n}}}/\eus{I}_{\rm GT} \simeq \C[\partial_x,\partial_y,z]$;
  \item[8)] $\eus{I}_{{\rm GT}^*}=(\partial_x,\partial_y,z)$, $\eus{A}^\mfrak{g}_{\widebar{\mfrak{n}}}/\eus{I}_{{\rm GT}^*} \simeq \C[x,y,\partial_z]$.
\end{enumerate}
Hence, we have the complete list of left ideals of the Weyl algebra $\eus{A}^\mfrak{g}_{\widebar{\mfrak{n}}}$ of the given type, where we consider either $x_\alpha$ or $\partial_{x_\alpha}$ for each $\alpha \in \Delta^+$.
\medskip

\theorem{\label{thm:irreducibility sl(3,C)}
\emph{Let $\lambda \in \mfrak{h}^*$. Then $W^\mfrak{g}_\mfrak{b}(\lambda)$ is a simple $\mfrak{g}$-module if and only if $\lambda(h_{\alpha_1}), \lambda(h_{\alpha_2}) \notin \N_0$ and $\lambda(h_\theta) \notin \Z$.}}

\proof{Let $M$ be a nonzero $\mfrak{g}$-submodule of $W^\mfrak{g}_\mfrak{b}(\lambda)$ for $\lambda \in \mfrak{h}^*$. Then $M$ is a weight $\mfrak{h}$-module, since $W^\mfrak{g}_\mfrak{b}(\lambda-\rho)$ is a weight $\mfrak{h}$-module. Moreover, we have $\pi_\lambda(f_\theta)=-\partial_z$, which implies that $M$ contains a nonzero weight vector $v$ satisfying $\partial_z v=0$. Hence, the weight vector $v$ is of the form $v=\partial_x^a \partial_y^b$ for some $a,b \in \N_0$. Further, by Proposition \ref{prop:pi action sl(2,C)-triple} and Proposition \ref{prop:pi action nilradical} we have
\begin{align*}
  \pi_\lambda(e_\theta)& = z(\pi_\lambda(h_\theta)-z\partial_z) + \sigma_\lambda(e_\theta)
\end{align*}
and
\begin{align*}
  \pi_\lambda(e_\alpha)=z\pi_\lambda([e_\alpha,f_\theta]) + \sigma_\lambda(e_\alpha)
\end{align*}
for $\alpha \in \Delta^+_\theta$. As we have $[\pi_\lambda(h_\theta) ,z]=2z$, we may write
\begin{align*}
  (\pi_\lambda(h_\theta)-2)\pi_\lambda(e_\alpha)= z\pi_\lambda(h_\theta)\pi_\lambda([e_\alpha,f_\theta]) + (\pi_\lambda(h_\theta)-2)\sigma_\lambda(e_\alpha)
\end{align*}
and
\begin{align*}
  \pi_\lambda(e_\theta)\pi_\lambda([e_\alpha,f_\theta])= z(\pi_\lambda(h_\theta)-z\partial_z)\pi_\lambda([e_\alpha,f_\theta]) + \sigma_\lambda(e_\theta)\pi_\lambda([e_\alpha,f_\theta])
\end{align*}
for $\alpha \in \Delta^+_\theta$. Since $[\pi_\lambda(f_\theta), \pi_\lambda([e_\alpha,f_\theta])]=0$ for $\alpha \in \Delta^+_\theta$ and $\partial_z v=0$, we get
\begin{align*}
  (\pi_\lambda(h_\theta)-2)\pi_\lambda(e_\alpha)v&= z\pi_\lambda(h_\theta)\pi_\lambda([e_\alpha,f_\theta])v + (\pi_\lambda(h_\theta)-2)\sigma_\lambda(e_\alpha)v, \\
  \pi_\lambda(e_\theta)\pi_\lambda([e_\alpha,f_\theta])v&= z\pi_\lambda(h_\theta)\pi_\lambda([e_\alpha,f_\theta])v + \sigma_\lambda(e_\theta)\pi_\lambda([e_\alpha,f_\theta])v,
\end{align*}
which implies that
\begin{align*}
  ((\pi_\lambda(h_\theta)-2)\pi_\lambda(e_\alpha)- \pi_\lambda(e_\theta)\pi_\lambda([e_\alpha,f_\theta]))v= ((\pi_\lambda(h_\theta)-2)\sigma_\lambda(e_\alpha)- \sigma_\lambda(e_\theta)\pi_\lambda([e_\alpha,f_\theta]))v
\end{align*}
for $\alpha \in \Delta^+_\theta$. If we denote
\begin{align*}
 w_1=&((\pi_\lambda(h_\theta)-2)\sigma_\lambda(e_{\alpha_1}) - \sigma_\lambda(e_\theta)\pi_\lambda([e_{\alpha_1},f_\theta]))v, \\
 w_2=&((\pi_\lambda(h_\theta)-2)\sigma_\lambda(e_{\alpha_2}) - \sigma_\lambda(e_\theta)\pi_\lambda([e_{\alpha_2},f_\theta]))v,
\end{align*}
then we have $w_1, w_2 \in M$ and $\pi_\lambda(f_\theta)w_1=0$, $\pi_\lambda(f_\theta)w_2=0$.
We may write
\begin{align*}
  w_1&= \big((\pi_\lambda(h_\theta)-2)\sigma_\lambda(e_{\alpha_1}) + \sigma_\lambda(e_\theta)\pi_\lambda(f_{\alpha_2})\big)v \\
  &= \big((x\partial_x+y\partial_y+\lambda_1+\lambda_2) x(x\partial_x - {\textstyle {1\over 2}} y\partial_y + \lambda_1+1) - {\textstyle {1 \over 2}}xy(x\partial_x-y\partial_y+\lambda_1-\lambda_2)\partial_y\big)v \\
  &=x(x\partial_x+\lambda_1+\lambda_2+1)(x\partial_x+\lambda_1+1)v = x(\partial_x x+\lambda_1+\lambda_2)(\partial_x x+\lambda_1)v
\end{align*}
and similarly
\begin{align*}
  w_2&=((\pi_\lambda(h_\theta)-2)\sigma_\lambda(e_{\alpha_2}) - \sigma_\lambda(e_\theta)\pi_\lambda(f_{\alpha_1}))v \\
  &= \big((x\partial_x+y\partial_y+\lambda_1+\lambda_2) y(y\partial_y - {\textstyle {1\over 2}} x\partial_x + \lambda_2+1) + {\textstyle {1 \over 2}}xy(x\partial_x-y\partial_y+\lambda_1-\lambda_2)\partial_x\big)v \\
  &=y(y\partial_y+\lambda_1+\lambda_2+1)(y\partial_y+\lambda_2+1)v= y(\partial_yy+\lambda_1+\lambda_2)(\partial_y y+\lambda_2)v,
\end{align*}
where we used Theorem \ref{thm:embedding sl(3,C)} and notation  $\lambda_1=\lambda(h_{\alpha_1})$ and $\lambda_2=\lambda(h_{\alpha_2})$. Since $v=\partial_x^a\partial_y^b$ for some $a,b \in \N_0$, we obtain
\begin{align*}
  w_1 = -a(\lambda_1+\lambda_2-a)(\lambda_1-a) \partial_x^{a-1} \partial_y^b \qquad \text{and} \qquad  w_2=  -b(\lambda_1+\lambda_2-b)(\lambda_2-b)\partial_x^a \partial_y^{b-1},
\end{align*}
which implies that $a\partial_x^{a-1} \partial_y^b, b\partial_x^a \partial_y^{b-1} \in M$ provided $\lambda_1, \lambda_2, \lambda_1+\lambda_2 \notin \N$. As a consequence, we get that $1 \in M$ if $\lambda_1,\lambda_2,\lambda_1+\lambda_2 \notin \N$. Since $W^\mfrak{g}_\mfrak{b}(\lambda)$ is generated by the vector $1 \in W^\mfrak{g}_\mfrak{b}(\lambda)$ for $\lambda_1 + \lambda_2 \notin -\N_0$, as follows from the proof of Theorem \ref{thm:cyclic}, we obtain that $W^\mfrak{g}_\mfrak{b}(\lambda)$ is a simple $\mfrak{g}$-module if $\lambda_1, \lambda_2 \notin \N$ and $\lambda_1+\lambda_2 \notin \Z$.

Now, let us assume that $\lambda_1+\lambda_2 \in -\N_0$. Let $M$ be the $\mfrak{g}$-submodule of $W^\mfrak{g}_\mfrak{b}(\lambda)$ generated by the vector $1 \in W^\mfrak{g}_\mfrak{b}(\lambda)$. For a vector $v \in M$ we may write $v=\pi_\lambda(u)1$, where $u \in U(\mfrak{g})$, which gives us $\pi_\lambda(e_\theta)v= \pi_\lambda([e_\theta,u])1 + \pi_\lambda(u)\pi_\lambda(e_\theta)1$. Since $\ad(e_\theta)$ is locally nilpotent on $U(\mfrak{g})$ and $\pi_\lambda(e_\theta)^{-\lambda_1-\lambda_2+1}1=0$, we obtain that $\pi_\lambda(e_\theta)$ is locally nilpotent on $M$. However, we have $\pi_\lambda(e_\theta)^k z^{-\lambda_1-\lambda_2+1}=k! z^{-\lambda_1-\lambda_2+k+1}$ for $k \in \N_0$, which implies that $z^{-\lambda_1-\lambda_2+1} \notin M$. Hence, we get that $W^\mfrak{g}_\mfrak{b}(\lambda)$ is not a simple $\mfrak{g}$-module.

Further, let us assume that $\lambda_1 \in \N$. Then there exists a homomorphism
\begin{align*}
  \varphi_\lambda \colon M^\mfrak{g}_\mfrak{b}(s_1(\lambda +\rho) - \rho) \rarr M^\mfrak{g}_\mfrak{b}(\lambda)
\end{align*}
of Verma modules. As the Verma module $M^\mfrak{g}_\mfrak{b}(\lambda)$ for $\lambda \in \mfrak{h}^*$ is realized as $\eus{A}^\mfrak{g}_{\widebar{\mfrak{n}}}/\eus{I}_{\rm V} \simeq \C[\partial_x,\partial_y,\partial_z]$, where the $\mfrak{g}$-module structure on $\eus{A}^\mfrak{g}_{\widebar{\mfrak{n}}}/\eus{I}_{\rm V}$ is given through the homomorphism $\pi_{\lambda+\rho} \colon \mfrak{g} \rarr \eus{A}^\mfrak{g}_{\widebar{\mfrak{n}}}$ of Lie algebras, we have
\begin{align*}
  \pi_\lambda(a)\varphi_\lambda(v)=\varphi_\lambda(\pi_{s_1\lambda}(a)v)
\end{align*}
for all $v \in \C[\partial_x,\partial_y,\partial_z]$ and $a \in \mfrak{g}$. Moreover, we have $\varphi_\lambda = (\partial_x + {1 \over 2} y\partial_z)^{\lambda_1}$, as follows immediately from \cite[Theorem 2.1]{Krizka-Somberg2015b}, which gives us $\pi_\lambda(a) \circ \varphi_\lambda = \varphi_\lambda \circ \pi_{s_1\lambda}(a)$ in the Weyl algebra \smash{$\eus{A}^\mfrak{g}_{\widebar{\mfrak{n}}}$} for all $a \in \mfrak{g}$. Therefore, we obtain a homomorphism
\begin{align*}
  \psi_\lambda \colon W^\mfrak{g}_\mfrak{b}(s_1 (\lambda + \rho)  - \rho) \rarr W^\mfrak{g}_\mfrak{b}(\lambda)
\end{align*}
of Gelfand--Tsetlin modules, such that $\psi_\lambda=\varphi_\lambda$ if we realize the Gelfand--Tsetlin module $W^\mfrak{g}_\mfrak{b}(\lambda)$ for $\lambda \in \mfrak{h}^*$ as $\eus{A}^\mfrak{g}_{\widebar{\mfrak{n}}}/\eus{I}_{\rm GT} \simeq \C[\partial_x,\partial_y,z]$.
As $\psi_\lambda$ is not surjective, we obtain that $W^\mfrak{g}_\mfrak{b}(\lambda)$ is not a simple $\mfrak{g}$-module. By a similar argument, we prove it in the case $\lambda_2 \in \N$ or $\lambda_1 + \lambda_2 \in \N$. If $\lambda_2 \in \N$, then there exists a homomorphism
\begin{align*}
  \varphi_\lambda \colon M^\mfrak{g}_\mfrak{b}(s_2(\lambda + \rho)  - \rho) \rarr M^\mfrak{g}_\mfrak{b}(\lambda)
\end{align*}
of Verma modules with $\varphi_\lambda = (\partial_y - {1 \over 2}y \partial_z)^{\lambda_2}$.   Finally, if $\lambda_1 + \lambda_2 \in \N$, then there exists a homomorphism
\begin{align*}
  \varphi_\lambda \colon M^\mfrak{g}_\mfrak{b}(s_1s_2s_1 (\lambda + \rho) - \rho) \rarr M^\mfrak{g}_\mfrak{b}(\lambda)
\end{align*}
of Verma modules with $\varphi_\lambda = \prod_{j=0}^{k-1} \big((\partial_x+{1 \over 2} y \partial_z)(\partial_y-{1\over 2}x \partial_z)-{1 \over 2}(\lambda_1-\lambda_2 -k+2j )\partial_z\big)$, where $k=-\lambda_1-\lambda_2$.

Therefore, the Gelfand--Tsetlin module $W^\mfrak{g}_\mfrak{b}(\lambda)$ is a simple $\mfrak{g}$-module if and only if $\lambda_1, \lambda_2 \notin \N$ and $\lambda_1+\lambda_2 \notin \Z$.}
\vspace{-2mm}


\begin{appendices}

\section{Generalized eigenspace decomposition}
\label{app:eigenspace decomposition}

For the reader's convenience, we summarize and prove several important facts concerning the generalized eigenspace decomposition.
\medskip

Let $V$ be a complex vector space and let $T \colon V \rarr V$ be a linear mapping. For each $\lambda \in \C$ we set
\begin{align}
  V_\lambda = \{v \in V;\, (\exists k \in \N)\, (T-\lambda \id)^kv=0\}
\end{align}
and call it the generalized eigenspace of $T$ with eigenvalue $\lambda$. When $V_\lambda \neq \{0\}$, we say that $\lambda$ is an eigenvalue of $T$ and the elements of $V_\lambda$ are called
generalized eigenvectors with eigenvalue $\lambda$. Further, let us note that
\begin{align}
  \sum_{\lambda \in \C} V_\lambda = \bigoplus_{\lambda \in \C} V_\lambda,
\end{align}
however it does not hold, in general, that $V = \bigoplus_{\lambda \in \C} V_\lambda$.
\medskip

Let $V = \bigoplus_{r\in \N_0}\! V_r$ be an $\N_0$-graded complex vector space. Then we denote by $p_r \colon V \rarr V_r$ for $r \in \N_0$ the canonical projection with respect to the direct sum decomposition. Furthermore, we define
\begin{align}
  F_k V = \bigoplus_{r=0}^k V_r
\end{align}
for all $k \in \N_0$, which gives us an increasing filtration $\{F_kV\}_{k \in \N_0}$ on $V$. \medskip

\lemma{\label{lem:eigenspace decomposition}
Let $V=\bigoplus_{r \in \N_0}\! V_r$ be an $\N_0$-graded complex vector space and let $\dim V_r < \infty$ for all $r \in \N_0$. Further, let $T \colon V \rarr V$ be a linear mapping satisfying
\begin{align}
  T(F_kV) \subset F_kV
\end{align}
for all $k \in \N_0$. Then we have a decomposition
\begin{align}
  V = \bigoplus_{\lambda \in \C} V_\lambda
\end{align}
called the generalized eigenspace decomposition of $V$ with respect to $T$.}

\proof{We already know that $\bigoplus_{\lambda \in \C} V_\lambda \subset V$. Let us assume that $v \in V$. Then there exists $k \in \N_0$ such that $v \in F_k V$. Since $T(F_kV) \subset F_kV$ and $\dim F_k V < \infty$, we have the generalized eigenspace decomposition of $F_kV$ with respect to $T$. Therefore, we obtain a decomposition $v = \sum_{\lambda \in \C} v_\lambda$, were $v_\lambda \in V_\lambda \cap F_kV$, which gives us the required statement.}

\theorem{\label{thm:generalized eigenspace isomorphism}
\emph{Let $V=\bigoplus_{r \in \N_0}\! V_r$ be an $\N_0$-graded complex vector space and let $\dim V_r < \infty$ for all $r \in \N_0$. Further, let $T_s, T_n \colon V \rarr V$ be linear mappings satisfying
\begin{align}
  T_s(V_r) \subset V_r \qquad \text{and} \qquad T_n(F_rV) \subset F_{r-1}V \label{eq:Ts and Tn condition}
\end{align}
for all $r \in \N_0$. Then there exists an isomorphism
\begin{align}
  \varphi \colon V \rarr V
\end{align}
of vector spaces such that
\begin{align}
  \varphi(V_\lambda^s) = V_\lambda^{s+n}
\end{align}
for all $\lambda \in \C$, where
\begin{align}
  V = \bigoplus_{\lambda \in \C} V_\lambda^s \qquad \text{and} \qquad V = \bigoplus_{\lambda \in \C} V_\lambda^{s+n}
\end{align}
are the generalized eigenspace decompositions of $V$ with respect to $T_s$ and $T_s+T_n$, respectively.}}

\proof{Since $T_s(F_kV) \subset F_kV$ and $(T_s+T_n)(F_kV) \subset F_kV$ for all $k \in \N_0$, as follows immediately from \eqref{eq:Ts and Tn condition}, we obtain by Lemma \ref{lem:eigenspace decomposition} the generalized eigenspace decompositions
\begin{align*}
  V = \bigoplus_{\lambda \in \C} V_\lambda^s \qquad \text{and} \qquad V = \bigoplus_{\lambda \in \C} V_\lambda^{s+n}
\end{align*}
of $V$ with respect to $T_s$ and $T_s+T_n$, respectively.

Now, we construct a linear mapping $\varphi \colon V \rarr V$. Since we have the direct sum decomposition
\begin{align*}
  V = \bigoplus_{r \in \N_0} \bigoplus_{\lambda \in \C}\, (V_r \cap V_\lambda^s),
\end{align*}
it is enough to define a linear mapping $\varphi$ on $V_r \cap V_\lambda^s$, which we shall denote by $\varphi_{r,\lambda}$, for all $r \in \N_0$, $\lambda \in \C$ and then extend linearly to $V$.
To define $\varphi_{r,\lambda}$, let $w \in V_r$ be a generalized eigenvector of $T_s$ with eigenvalue $\lambda \in \C$. As $w \in F_r V$ and $(T_s+T_n)(F_rV) \subset F_rV$, we get a decomposition
\begin{align*}
w = \sum_{\mu \in \C} w_\mu^{s+n},
\end{align*}
were $w_\mu^{s+n} \in F_rV \cap V^{s+n}_\mu$ for all $\mu\in \C$. Applying the projection $p_r \colon V \rarr V_r$ on the previous equation, we obtain
\begin{align*}
  w = \sum_{\mu \in \C} p_r(w_\mu^{s+n}),
\end{align*}
since $p_r(w)=w$. As $w_\mu^{s+n} \in F_rV \cap V_\mu^{s+n}$, there exists $k \in \N$ such that $(T_s+T_n-\mu)^kw_\mu^{s+n}=0$ and therefore we may write
\begin{align*}
 0= p_r((T_s+T_n-\mu)^kw_\mu^{s+n})= p_r((T_s-\mu)^kw_\mu^{s+n})=(T_s-\lambda)^kp_r(w_\mu^{s+n}),
\end{align*}
where we used \eqref{eq:Ts and Tn condition}, which gives us that \smash{$p_r(w_\mu^{s+n})$} is a generalized eigenvector of $T_s$ with eigenvalue $\mu$ for all $\mu \in \C$. Hence, the decomposition \smash{$w = \sum_{\mu \in \C} p_r(w_\mu^{s+n})$} implies that $p_r(w_\mu^{s+n})=0$ for $\mu \neq \lambda$ and \smash{$p_r(w_\lambda^{s+n})=w$}. Therefore, we may set
\begin{align*}
  \varphi_{r,\lambda}(w) =w_\lambda^{s+n}
\end{align*}
for all $w \in V_r \cap V_\lambda^s$. The previous construction ensures that $\varphi_{r,\lambda} \colon V_r \cap V_\lambda^s \rarr F_r V \cap V_\lambda^{s+n}$ is a linear mapping.  Moreover, we have $\varphi(V_\lambda^s) \subset V_\lambda^{s+n}$ for all $\lambda \in \C$.

Further, since $\varphi(F_r V) \subset F_r V$ for all $r \in \N_0$, we obtain a linear mapping
\begin{align*}
  \gr_r\! \varphi \colon F_r V /F_{r-1}V \rarr F_r V / F_{r-1}V
\end{align*}
for all $r \in \N_0$. As $\gr_r\! \varphi = \id_{F_r V / F_{r-1}V}$ by construction, we get that $\varphi \colon V \rarr V$ is an isomorphism of vector spaces.}

\end{appendices}

\vspace{-2mm}


\section*{Acknowledgments}

V.\,F.\ is supported in part by CNPq (304467/2017-0) and by
Fapesp (2014/09310-5); L.\,K.\ is supported by Capes (88887.137839/2017-00).



\providecommand{\bysame}{\leavevmode\hbox to3em{\hrulefill}\thinspace}
\providecommand{\MR}{\relax\ifhmode\unskip\space\fi MR }
\providecommand{\MRhref}[2]{%
  \href{http://www.ams.org/mathscinet-getitem?mr=#1}{#2}
}
\providecommand{\href}[2]{#2}

\end{document}

Let us introduce a root basis of the niradical $\mfrak{n}$ by
\begin{align*}
  e_{\alpha_1}=
  \begin{pmatrix}
    0 & 1 & 0 \\
    0 & 0 & 0 \\
    0 & 0 & 0
  \end{pmatrix}\!, \qquad \qquad
  e_\theta&=
  \begin{pmatrix}
    0 & 0 & 1  \\
    0 & 0 & 0 \\
    0 & 0 & 0
  \end{pmatrix}\!, \qquad \qquad
  e_{\alpha_2}=
  \begin{pmatrix}
    0 & 0 & 0  \\
    0 & 0 & 1 \\
    0 & 0 & 0
  \end{pmatrix}
\end{align*}
and a root basis of the opposite nilradical $\widebar{\mfrak{n}}$ through
\begin{align*}
  f_{\alpha_1}=
  \begin{pmatrix}
    0 & 0 & 0 \\
    1 & 0 & 0 \\
    0 & 0 & 0
  \end{pmatrix}\!, \qquad \qquad
  f_\theta=
  \begin{pmatrix}
    0 & 0 & 0  \\
    0 & 0 & 0 \\
    1 & 0 & 0
  \end{pmatrix}\!, \qquad \qquad
  f_{\alpha_2}=
  \begin{pmatrix}
    0 & 0 & 0  \\
    0 & 0 & 0 \\
    0 & 1 & 0
  \end{pmatrix}\!.
\end{align*}
The coroots corresponding to the positive roots are given by
\begin{align*}
  h_{\alpha_1}=
  \begin{pmatrix}
    1 & 0 & 0 \\
    0 & -1 & 0 \\
    0 & 0 & 0
  \end{pmatrix}\!, \qquad \qquad
  h_\theta=
  \begin{pmatrix}
    1 & 0 & 0 \\
    0 & 0 & 0 \\
    0 & 0 & -1
  \end{pmatrix}\!, \qquad \qquad
  h_{\alpha_2}=
  \begin{pmatrix}
    0 & 0 & 0  \\
    0 & 1 & 0 \\
    0 & 0 & -1
  \end{pmatrix}\!.
\end{align*}

According to the classification of all simple Gelfand--Tsetlin modules for $\mfrak{sl}(3,\C)$ given in \cite{Futorny-Grantcharov-Ramirez2014}, we can identify $W^\mfrak{g}_\mfrak{b}(\lambda)$ with  ????? provided $\lambda(h_{\alpha_1}), \lambda(h_{\alpha_2}) \notin \N_0$ and $\lambda(h_\theta) \notin \Z$.